\documentstyle[12pt]{article}
\newcommand{\np}{}
\newcommand{\rbox}[1]{{\rm{#1}}}        
\newcommand{\bbox}[1]{{\bf{#1}}}
\newcommand{\bsym}[1]{{\bf{#1}}}
\newcommand{\rmt}{\rm}
\newcommand{\itt}{\it}
\newcommand{\bft}{\bf}
\expandafter\chardef\csname pre amsym.def at\endcsname=\the\catcode`\@
\catcode`\@=11
\def\newsymbol#1#2#3#4#5{\let\next@\relax
 \ifnum#2=\@ne\let\next@\msafam@\else
 \ifnum#2=\tw@\let\next@\msbfam@\fi\fi
 \mathchardef#1="#3\next@#4#5}
\def\hexnumber@#1{\ifcase#1 0\or 1\or 2\or 3\or 4\or 5\or 6\or 7\or 8\or
 9\or A\or B\or C\or D\or E\or F\fi}
\newfam\msafam
\newfam\msbfam
\edef\msafam@{\hexnumber@\msafam}
\edef\msbfam@{\hexnumber@\msbfam}
\def\Bbb#1{\fam\msbfam\relax#1}
\catcode`\@=\csname pre amsym.def at\endcsname
\newsymbol\subsetneqq 2324
\newsymbol\restriction 1316
\font\twlmsa=msam10 scaled \magstep1
\font\egtmsa=msam8
\def\xiimsa{\textfont\msafam=\twlmsa\scriptfont\msafam=\egtmsa} 
\font\twlmsb=msbm10 scaled \magstep1
\font\tenmsb=msbm10 
\def\xiimsb{\textfont\msbfam = \twlmsb}
\def\xmsb  {\textfont\msbfam = \tenmsb}
\newcommand{\xiieufm}{\font\teufm=eufm10 scaled \magstep1}
\newcommand{\xeufm}  {\font\teufm=eufm10}
\newcommand{\frak}[1]{\hbox{\teufm{#1}}}

\newcommand{\gM}  {{\bbb\frak{M}\bbb}}
\newcommand{\gN}  {{\bbb\frak{N}\bbb}}
\newcommand{\xiieusm}{\font\teusm=eusm10 scaled \magstep1
\font\seusm=eusm9 \font\zeusm=eusm7}
\newcommand{\xeusm}{\font\teusm=eusm10 \font\seusm=eusm8 
\font\zeusm=eusm6}
\newcommand{\skr}[1]{{\mathchoice {\hbox{\teusm{#1}}} 
{\hbox{\teusm{#1}}}{\hbox{\seusm{#1}}}{\hbox{\zeusm{#1}}} }}
\newcommand{\xiieurm}{\font\teurm=eurm10 scaled \magstep1
\font\seurm=eurm9 \font\zeurm=eurm7}
\newcommand{\xeurm}{\font\teurm=eurm10 \font\seurm=eurm8 
\font\zeurm=eurm6}
\newcommand{\cur}[1]{{\mathchoice {\hbox{\teurm{#1}}} 
{\hbox{\teurm{#1}}}{\hbox{\seurm{#1}}}{\hbox{\zeurm{#1}}} }}
\newcommand{\xiisf}{\font\tsf=cmss10 scaled \magstep1
\font\ssf=cmss9 \font\zsf=cmss9}
\newcommand{\xisf}{\font\tsf=cmss10 scaled \magstephalf 
\font\ssf=cmss8 \font\zsf=cmss8}
\newcommand{\xsf}{\font\tsf=cmss10 \font\ssf=cmss8 
\font\zsf=cmss8}
\newcommand{\erlf}[1]{{\mathchoice {\hbox{\tsf{#1}}} 
{\hbox{\tsf{#1}}}{\hbox{\ssf{#1}}}{\hbox{\zsf{#1}}} }}
\def\addto#1#2{
\ifx\zzone\undefined\let\zzone=#1\def#1{\zzone#2}\else
\ifx\zztwo\undefined\let\zztwo=#1\def#1{\zztwo#2}\else
\ifx\zzthr\undefined\let\zzthr=#1\def#1{\zzthr#2}\else
\fi\fi\fi
}
\addto\normalsize   {\xiimsb\xiimsa\xiieufm\xiieusm\xiieurm\xiisf}
\addto\footnotesize {\xeusm\xmsb\xeurm\xeufm\xsf}   
\addto\small        {\xisf}     
\newtheorem{theorem}{Theorem}
\newtheorem{assertion}[theorem]{Assertion}
\newtheorem{corollary}[theorem]{Corollary}
\newtheorem{definition}[theorem]{Definition}
\newtheorem{lemma}[theorem]{Lemma}
\newtheorem{proposition}[theorem]{Proposition}
\newtheorem{remark}[theorem]{Remark}
\newcommand{\TF}{\itt}
\newcommand{\bass}{\begin{assertion}\TF\ } 
\newcommand{\eass}{\end{assertion}}
\newcommand{\bcor}{\begin{corollary}\TF\ }
\newcommand{\ecor}{\end{corollary}}
\newcommand{\bdf} {\begin{definition}\rmt\ }
\newcommand{\edf} {\end{definition}} 
\newcommand{\ble} {\begin{lemma}\TF\ }
\newcommand{\ele} {\end{lemma}}
\newcommand{\bpro}{\begin{proposition}\TF\ } 
\newcommand{\epro}{\end{proposition}} 
\newcommand{\brem}{\begin{remark}\rmt\ }
\newcommand{\erem}{\end{remark}} 
\newcommand{\bte} {\begin{theorem}\TF\ }
\newcommand{\ete} {\end{theorem}}
\newcommand{\proof}{\noi{\bft Proof\hspace{1mm} }}
\newcommand{\qed}  {\hfill$\msur\Box\msur$} 
\newcommand{\qedd} {\hfill$\msur\dashv\msur$} 
\newcommand{\bay}{\begin{array}}
\newcommand{\eay}{\end{array}} 
\newcommand{\bce}{\begin{center}}
\newcommand{\ece}{\end{center}} 
\newcommand{\bde}{\begin{description}}
\newcommand{\ede}{\end{description}}
\newcommand{\ben}{\begin{enumerate}}
\newcommand{\een}{\end{enumerate}}
\newcommand{\bit}{\begin{itemize}}
\newcommand{\eit}{\end{itemize}}

\newcommand{\ZFC} {\bbox{ZFC}}

\newcommand{\bS}  {\bbox{S}}
\newcommand{\namex}  {\bbox{x}}
\newcommand{\namef}  {\bbox{f}}
\newcommand{\card}{\rbox{card}\,}
\newcommand{\dom} {\rbox{dom}\,}
\newcommand{\hc}  {\rbox{HC}}
\newcommand{\rL}  {\rbox{L}}
\newcommand{\ord} {\rbox{Ord}}
\newcommand{\Ord} {\ord}
\newcommand{\od}  {\rbox{OD}}
\newcommand{\ran} {\rbox{ran}\,}
\newcommand{\rod} {\hbox{\rmt R\hspace{1pt}-\hspace{1pt}OD}}
\newcommand{\rV}  {\rbox{V}}
\newcommand{\Vp}  {{\rbox{V}}^\ast}
\newcommand{\sma} [1]{#1\hbox{-}\rbox{\tt SM}}
\newcommand{\clps}[1]{#1\hbox{-}\rbox{\tt Clps}}
\newcommand{\ima} {{\hbox{\hspace{1pt}\rmt ''}}}

\newcommand{\abs} [1]{|\hspace{1pt}{#1}\hspace{1pt}|}
\newcommand{\poset}{p.\hspace{0.3em}o.\hspace{0.3em}set}
\newcommand{\al} {\alpha} 
\newcommand{\ba} {\beta} 
\newcommand{\ga} {\gamma} 
\newcommand{\da} {\delta}
\newcommand{\Da} {\Delta}
\newcommand{\kpa}{\kappa}
\newcommand{\La} {\Lambda} 
\newcommand{\la} {\lambda} 
\newcommand{\sg} {\sigma} 
\newcommand{\Sg} {\Sigma}
 
\newcommand{\vpi}{\varphi} 
\newcommand{\om} {\omega} 
\newcommand{\omi}{{\om_1}}
\newcommand{\Om} {\Omega}
\newcommand{\lom}{{<\om}} 

\newcommand{\fs}[2]{{\bsym\Sigma}^{#1}_{#2}}
\newcommand{\fp}[2]{{\bsym\Pi}^{#1}_{#2}}
\newcommand{\fd}[2]{{\bsym\Delta}^{#1}_{#2}}

\newcommand{\iSigma}{{\mathchar"7106}}
\newcommand{\iPi}   {{\mathchar"7105}}
\newcommand{\iDelta}{{\mathchar"7101}}
\newcommand{\is}[2]{\iSigma^{#1}_{#2}}
\newcommand{\ip}[2]{\iPi^{#1}_{#2}}
\newcommand{\id}[2]{\iDelta^{#1}_{#2}}
\newcommand{\ish}[1]{\is{\rbox{HC}}{#1}}
\newcommand{\iph}[1]{\ip{\rbox{HC}}{#1}}
\newcommand{\idh}[1]{\id{\rbox{HC}}{#1}}
\newcommand{\fsh}[1]{\fs{\rbox{HC}}{#1}}

\newcommand{\fdh}[1]{\fd{\rbox{HC}}{#1}}
\newcommand{\freepi}{{\char"19}}
\newcommand{\freesg}{{\char"1B}}
\newcommand{\bbb}{\hspace{0.5pt}} 
\newcommand{\dvoj}[1]{{\bbb{\Bbb #1}\bbb}}

\newcommand{\dF} {{\dvoj F}}

\newcommand{\dT} {{\dvoj T}}

\newcommand{\dX} {{\dvoj X}}

\newcommand{\cursp}{\hspace{0.5pt}}
\newcommand{\curi}[1]{{\cursp \cur{#1}\cursp}}
\newcommand{\ww}{\fsg}    
\newcommand{\zC}[1]{{\cur C}_{#1}}
\newcommand{\zP}   {{\cur P}}
\newcommand{\zX}[1]{{\cur X}_{#1}}
\newcommand{\fpi}{{\curi\freepi}}
\newcommand{\fsg}{{\curi\freesg}}
\newcommand{\zt}    {{t}}

\newcommand{\skrsp}{\hspace{0.5pt}}
\newcommand{\skri}[1]{{\skrsp \skr{#1}\skrsp}}
\newcommand{\cD} {{\skri D}}   
\newcommand{\cG} {{\cal G}}  
      
\newcommand{\cC}[1] {{\cal P}_{#1}}

\newcommand{\cN} {{\skri N}}
\newcommand{\cO} {{\skri O}}
\newcommand{\cP} {{\skri P}}
\newcommand{\oP} {{\skri P}^{\hbox{\tiny\rm OD}}}
\newcommand{\cW} {{\skri W}}
\newcommand{\cWt}{{\cW_{(2)}}}
\newcommand{\cX} {{\skri X}}
\newcommand{\emps}{\emptyset}
\newcommand{\sq}  {\subseteq}

\newcommand{\sneq}{\subsetneqq}
\newcommand{\cj}  {\;\,\&\;\,}
\newcommand{\orr} {\,\;{\textstyle\bigvee}\;}
\newcommand{\lra} {\longrightarrow} 
\newcommand{\llra}{\longleftrightarrow} 
 
\newcommand{\res} {{\mathbin{\hspace{0.15ex}\restriction\hspace{0.15ex}}}}
\newcommand{\we}  {{\mathbin{\hspace{0pt}^\wedge}}}

\newcommand{\<}{\leq}
\newcommand{\ti}{\times}
\newcommand{\dm}{$$}
\newcommand{\ang} [1]{\langle #1\rangle}
\newcommand{\ans} [1]{\{\hspace{0.2mm}#1\hspace{0.2mm}\}}

\newcommand{\dd}[2]{\hpsur\hbox{\mathsurround=0mm${#1}$-#2}}
\newcommand{\itla}[1]{\item\label{#1}}
\newcommand{\hoc}{\mathbin{\<_{\rmt c}}}
\newcommand{\emc}{\mathbin{\sqsubseteq_{\rmt c}}}
\newcommand{\hocn}{\mathbin{\<_{\rmt cn}}}
\newcommand{\emcn}{\mathbin{\sqsubseteq_{\rmt cn}}}
\newcommand{\sfbox}[1]{{\erlf #1}}   
\newcommand{\rG} {\mathbin{\sfbox G}}
\newcommand{\rGw}{\rG}    
\newcommand{\rC} {\mathbin{\sfbox C}}
\newcommand{\rK} {\mathbin{\sfbox K}}
\newcommand{\Go}   {\mathbin{{\cG}_0}} 
\newcommand{\Gph}  {\mathbin{{\sfbox G}_\Phi}}

\newcommand{\Gp}   {\mathbin{{\sfbox G}^\ast}}
\newcommand{\rQ}[1]{\mathbin{{\sfbox Q}_{#1}}}

\newcommand{\nG} {\mathbin{\not{\hspace{-0.25em}\sfbox G}}}

\newcommand{\spp}{\mathbin{\widehat{\hspace{1ex}}}}

\newcommand{\doS}{{\bS}}   
\newcommand{\col}[1]{{#1}^\lom}   
\newcommand{\smf}[1]{{\cal P}_{#1}}
\newcommand{\smfu}[1]{{\cal P}^{#1}}
\newcommand{\smfl}[1]{{\cal P}_{<#1}}
\newcommand{\smfe}[1]{{\cal P}_{\<#1}}
\newcommand{\etal}{{\it et al\/}.}
\newcommand{\etc}{{\it etc\/}.}
\newcommand{\eg}{{\it e.\hspace{0.4ex}g\/}.}
\newcommand{\ie}{{\it i.\hspace{0.4ex}e\/}.}
\newcommand{\vs}{{\it vs\/}.}
\newcommand{\gfu}[1]{{\dF}_{#1}}      
\newcommand{\trb}[1]{{\dT}_{#1}}
\newcommand{\msur} {\hspace{-1\mathsurround}}

\newcommand{\hpsur}{\hspace{0.5\mathsurround}}
\newcommand{\psur} {\hspace{\mathsurround}}
\newcommand{\noi}  {\noindent}
\newcommand{\vom}  {\vspace{1mm}}
\newcommand{\hspep}{\vspace{3mm}}
\newcommand{\its}  {\vspace{-1mm}}
\begin{document}

\normalsize

\title{On a dichotomy related to colourings of definable graphs 
in generic models}

\author{Vladimir Kanovei
\thanks{Moscow Transport Engineering Institute}
\thanks{{\tt kanovei@math.uni-wuppertal.de} \ and \ 
{\tt kanovei@mech.math.msu.su}
}
\thanks{The author acknowledges the support from AMS, DFG,  
NWO, and universities of Wuppertal and Amsterdam.} 
}
\date{May 1996} 
\maketitle
\normalsize

\vfill

\begin{abstract}\vspace{2mm} 
\noi
\parskip=1mm
We prove that in the Solovay model every $\od$ graph $\rG$ 
on reals satisfies one and only one of the following two 
conditions: $(\hbox{I})\msur$ $\rG$ admits an $\od$ colouring by 
ordinals; $(\hbox{II})\msur$ there exists a continuous 
homomorphism of $\Go$ into $\rG,$ where $\Go$ is a certain 
${\bbox{F}}_\sg$ locally countable graph which is not $\rod$ 
colourable by ordinals in the Solovay model. If the graph $\rG$ is 
locally countable or acyclic then $(\hbox{II})$ can be 
strengthened by the requirement that the homomorphism is 
a $1-1$ map, \ie\ an embedding.

As the second main result we prove that $\fs12$ graphs admit the 
dichotomy $(\hbox{I})$ \vs\ $(\hbox{II})$ in set--generic 
extensions of the constructible universe $\rL$ (although now 
$(\hbox{I})$ and  $(\hbox{II})$ may be in general compatible). 
In this case $(\hbox{I})$ 
can be strengthened to the existence of a $\fd13$ colouring by 
countable ordinals provided the graph is locally countable. 

The proofs are based on a topology generated by $\od$ sets.

\end{abstract}

\vfill
$\,$

\newpage

\subsection*{Introduction}

A new direction in the classical domain of graph colouring 
was discovered by Kechris, Solecki, and Todorcevic~\cite{kst}. 
They found that graphs on reals drastically 
change their behaviour in the case when a definable 
colouring is requested. For instance there exists a 
${\bbox{F}}_\sg$ graph $\Go$ on reals (rather a family of 
graphs generated by a common method), acyclic (therefore 
colourable by only two colours) and locally countable, 
which does not admit even a countable {\it Borel\/} colouring. 

Moreover, it is proved in \cite{kst} that $\Go$ is in a 
certain way minimal among all analytic graphs which are not 
countably Borel colourable. 

To formulate this result consistently, let us recall some 
notation. We refer to \cite{kst} for a more substantial 
review.

A {\it graph\/} on a set $X$ (typically $X$ is a set of reals) 
is any set $\rG\sq X^2$ such that 
$x\nG x$ and $x\rG y\,\llra\,y\rG x$ for all $x,\,y\in X$. 

A graph $\rG$ is {\it acyclic\/} iff there exists no chain of the 
form \/ $x_0\rG x_1\rG x_2 \rG ... \rG x_n$ with\/ $x_0=x_n$ 
and\/ $n\geq 3.$ A graph $\rG$ is {\it locally countable\/} iff  
every vertex has at most countably many neighbours in $\rG$. 

A {\it colouring\/} of a graph $\rG$ on $X$ is a function $c$ 
defined on $X$ so that $x\rG y$ implies $c(x)\not=c(y).$ 
A colouring $c$ is {\it countable\/} iff $c(x)\in\om$ for all 
$x.$ In this case $c$ is a {\it Borel\/} colouring iff every 
pre--image $c^{-1}(n)$ $\msur (n\in\om)$ is a Borel set.
We shall also consider {\it colourings by ordinals\/}, 
meaning that $c(x)\in\Ord$ for all $x$. 

Let $\rG$ and $\rG'$ be graphs on topological spaces $X$ and 
$X'$ respectively. A {\it homomorphism\/} of $\rG'$ into $\rG$ 
is a continuous function $H:X'\,\lra\,X$ such that 
$x\rG' y\,\lra H(x)\rG H(y)$ for all $x,\,y\in X.$ An 
{\it embedding\/} is an injective homeomorphism, that is, 
a continuous $1-1$ function with the same property (but 
$x\rG' y\,\longleftarrow\,H(x)\rG H(y)$ is not requested). 

One writes ${\rG'\hoc\rG}$ if there exists a homeomorphism of 
$\rG'$ into $\rG.$ One writes ${\rG'\emc\rG}$ if there exists an 
embedding of $\rG'$ into $\rG.$ ($\hspace{0mm}_{\rmt c}$ stands 
for ``continuous''.) Even ${\rG'\hoc\rG}$ suffices for 
$\rG$ to inherit uncolourability properties from~$\rG'$. 

Let us fix a sequence $\ang{\ww_n:n\in\om}$ of 
binary tuples $\ww_n\in 2^n$ such that for any $s\in 2^\lom$ 
there exists $n$ such that $s\sq\ww_n$ (the {\it density\/}).
We define $\Go,$ a graph on the {\it Cantor space\/} $\cD=2^\om,$ 
as follows: $x\Go y$ iff there exists $n$ such that 
\dm
x\res n=y\res n=\ww_n\,,\hspace{4mm}x(n)=1-y(n)\,,
\hspace{4mm}\hbox{and}\hspace{4mm}
x(k)=y(k)\hspace{3mm}\hbox{for all}\hspace{3mm}k>n\,.~\footnote
{\ The essential here properties of $\Go$ do not depend on the 
particular choice of the tuples $\ww_n.$ However it seems unknown 
whether different graphs $\Go$ are isomorphic.} 
\dm
One easily proves that $\Go$ is a locally countable acyclic 
${\bbox{F}}_\sg$ graph on $\cD$ which does not 
admit a countable Borel colouring. (We refer to~\cite{kst}, see 
also Subsection~\ref{incomp}; the density is essential only for 
the Borel uncolourability.) Therefore any graph $\rG$ such that 
$\Go\hoc \rG$ does not admit a countable Borel colouring as well. 
The following theorem shows that this is also a necessary condition.

\bte 
\label{tkst}
{\rmt (Theorem 6.3 in \cite{kst})} \ 
Every analytic graph\/ $\rG$ on reals 
satisfies one and only one of the  
two statements,\/ {\rmt\ref{1111}} and\/ {\rmt\ref{2222}}$:$ 
\ben
\itemsep=1mm
\def\theenumi{$\langle{\hbox{\Roman{enumi}}}\rangle$}
\def\labelenumi{{\rmt\theenumi}}
\itla{1111}\msur 
$\rG$ admits a countable Borel colouring$.$

\itla{2222}\msur $\Go\hoc \rG.$ Moreover,\/ $\Go\emc\rG$ 
in each of the two cases$:$ 
$a)\msur$ $\rG$ is locally countable$;$ 
$b)\msur$ $\rG$ is acyclic.  
\een
\ete
\noi
It is quite usual in descriptive set theory that theorems 
proved for objects of initial projective classes (\eg\  
for Borel sets) generalize in the Solovay model~\footnote
{\rmt\ By {\it Solovay model\/} we mean a generic extension 
$\rL[G]$ of the constructible universe $\rL$ via a certain 
notion of forcing $\smfu\la,$ as defined by 
Solovay~\cite{sol}. See Section~\ref{somo} 
for more information.} 
on all $\rod$ (real--ordinal--definable) sets of reals. 
The following theorem is such a generalization of 
Theorem~\ref{tkst}.

\bte
\label{main}
In the Solovay model, every\/ $\rod$ graph\/ $\rG$ on reals 
satisfies one and only one of the  
two statements,\/ {\rmt\ref1} and\/ {\rmt\ref2}$:$ 
\ben
\itemsep=1mm
\def\theenumi{\hskip2pt{(\Roman{enumi})}\hskip2pt}
\def\labelenumi{{\rmt\theenumi}}
\itla{1}\msur 
$\rG$ admits a\/ $\rod$ colouring by countable ordinals. 

\itla{2}\msur $\Go\hoc \rG.$ Moreover,\/ $\Go\emc\rG$ 
in each of the two cases$:$ 
$a)\msur$ $\rG$ is locally countable$;$ 
$b)\msur$ $\rG$ is acyclic.  
\een
\ete
\noi
(It is pointed out in \cite{kst} with reference to Woodin that a 
dichotomy of this type is a consequence of some determinacy 
hypothesis.) 

Two open problems should be mentioned. 

{\it Problem 1\/}. Characterize the case when $\rG$ is 
{\it countably\/} $\rod$ colourable.

{\it Problem 2\/}. Prove Theorem~\ref{main} with \ref{2} in 
the form ${\Go\emc\rG}$ unconditionally. (This is also an 
open problem for Theorem~\ref{tkst}.)

A partial result in concern of the latter problem will be 
proved in Section~\ref{3rd}. Let a {\it curve\/} be a binary 
relation $\rK$ satisfying at least one of the two following 
requirements: $x\rK y\cj x\rK y'\,\lra\,y=y'$ or 
$x\rK y\cj x'\rK y\,\lra\,x=x'.$ (In other words $\rK$ is the 
graph of a partial function, either $y=f(x)$ or $x=f(y)$.) 
Let a {\it special\/} graph mean: a union of the form 
$\rG=\bigcup_{\xi<\omi}\rK_\xi,$ where each $\rK_\xi$ is 
a curve and the map $\xi\,\longmapsto\,\rK_\xi$ is $\rod$.

\bte
\label{spe}
To prove Theorem~\ref{main} for all\/ $\rod$ graphs with\/ 
{\rmt\ref{2}} in the unconditional form\/ ${\Go\emc\rG},$ it 
suffices to obtain this result for special graphs. 
\ete
This reduction corresponds to a result for Borel graphs 
(namely, reduction to the case when $\rG$ is  
a union of countably many Borel curves), absent in 
\cite{kst} but discussed between A.~S.~Kechris and the 
author.

It is also quite usual that theorems on Borel sets have 
(suitably weaker) counterparts for sets of initial 
projective levels, true in $\ZFC$ or in a quite reasonable 
hypothesis of broad nature, as \eg\ sharps or 
determinacy hypotheses. It has been recently found that 
the assumption that the universe $\rV$ of all sets is a 
set--generic extension of the constructible universe $\rL,$ 
helps to prove theorems on analytic ($\fs11$) equivalence 
relations still unproved in $\ZFC$ alone. (We refer to a 
theorem on thin $\fs11$ relations of Hjorth~\cite{h-thin} 
and a dichotomy theorem of Glimm -- Effros type proved for 
$\fs11$ relations by Kanovei~\cite{k-s11}. It is remarkable 
that these results have also been proved in $\ZFC$ plus the 
sharps hypothesis (Hjorth~\cite{h-thin}, Hjorth and 
Kechris~\cite{hk}) or the existence of a weakly compact 
cardinal (Friedman~\cite{sdf}), but a proof in $\ZFC$ alone 
is so far unknown.) The following theorem contributes in 
this type of results.

\bte
\label{tlast}
Suppose that the universe is a set--generic extension of 
the constructible universe. Then every\/ $\fs12$ graph\/ $\rG$ 
on reals satisfies one~\footnote
{\rmt\ In this case \ref{1v} and \ref{2v} are generally 
speaking compatible.}
of the two following statements$:$ 
\ben
\itemsep=1mm
\def\theenumi{\hskip2pt{(\Roman{enumi}{\mathsurround0mm$'$})}\hskip2pt}
\def\labelenumi{{\rmt\theenumi}}
\itla{1v}\msur 
$\rG$ admits a\/ $\rod$ colouring by ordinals. If in 
addition\/ $\rG$ is locally countable then\/ $\rG$ 
admits a\/ $\fdh 2$ colouring by countable ordinals.~\footnote
{\label{classes}\rmt\ $\hc$ is the set of all hereditarily 
countable sets. The ``boldface upright'' $\fsh n$ means 
\dd{\Sigma_n}definability in $\hc$ with reals being admitted 
as parameters. The ``lightface slanted'' $\ish n$ will  
mean parameter--free \dd{\Sigma_n}definability in $\hc.$ 
\dd\Pi notions are introduced in the same way, while $\Da$ 
in any form means simultaneously $\Sigma$ and $\Pi$.}

\itla{2v} --- exactly as\/ {\rmt\ref2} of Theorem~\ref{main}.
\een
\ete
Since \ref{1v} does not seem in any way adequate in the general 
case, it is an important problem to obtain something better than 
simply $\od$ in \ref{1v}. One more problem is to obtain matching 
counterexamples among $\fp12$ graphs. But the most interesting 
would be to prove the theorem in $\ZFC$ (not assuming that the 
universe is a set--generic extension of $\rL$). 

The exposition starts with a review of some facts in the theory of 
the Solovay and more elementary collapse models in sections \ref{somo} 
and \ref{reals}. We also introduce a topology generated 
by $\od$ sets (a counterpart of the Gandy -- Harrington topology) 
and study an associated forcing. 

Section \ref{1st} reveals the nature of the dichotomy. 
Concentrating on $\od$ graphs, we 
consider the case when the given $\od$ graph $\rG$  
satisfies the property that every real belongs to an $\od$ 
\dd\rG discrete set $X,$ and prove that this leads to 
\ref{1} in the Solovay model. (A set $X$ is \dd\rG{\it discrete\/} 
iff $\rG$ does not intersect $X\ti X$.)

Sections \ref{2nd} and \ref{3rd} handle the case when \underline{not} 
every real $x$ belongs to an $\od$ \dd\rG discrete set. We prove that 
then $\Go\hoc\rG,$ and moreover ${\Go\emc\rG}$ provided $\rG$ is 
locally countable or acyclic. This part of the 
proof of Theorem~\ref{main} utilizes an important splitting 
construction of Kechris \etal~\cite{kst}, but technically 
realized so that the Choquet games (as in \cite{kst}) are not involved. 

To explain the ideas behind the proof of Theorem~\ref{tlast}, let 
$\kpa$ be the \dd\rL cardinal of the notion of forcing 
$\zP\in\rL$ such that the universe $\rV$ is a \dd\zP generic 
extension of $\rL.$ Let $\rG$ be in fact a  
$\is12$ graph. Assume, for the sake of simplicity, that 
$\rV$ still contains an inaccessible cardinal $\Om>\kpa.$ 

Then $\rV$ can be embedded into a vider universe, $\Vp,$ which is a 
Solovay model over $\rL,$ so that $\Om=\omi^{\Vp}.$ (In fact a 
suitable collapse generic extension will be taken as $\Vp$.) The 
dichotomy \ref1 \vs\ \ref2 holds in $\Vp$ by Theorem~\ref{main}. 

If \ref2 holds in $\Vp$ then we obtain \ref{2v} in $\rV$ by an 
absoluteness argument based on the Shoenfield theorem. If \ref1 
holds in $\Vp$ then we get \ref{1v} in $\rV$ using the fact that 
\ref1 appears in $\Vp$ in the form of a statement somehow 
relativized to $\rL.$ 
This reasoning is presented in Section~\ref{s12}.

\subsubsection*{Acknowledgements} 

The author is pleased to thank 
A.\ S.\ Kechris for useful discussions and interesting 
information on the topic of this paper.

\subsubsection*{Notation}

We shall use the {\em Cantor set} $\cD=2^\om$ 
rather than the Baire space $\cN=\om^\om$ as the principal space 
in this paper. Elements of $\cD$ will be called {\it reals\/} 
below. 

In the rest we shall follow the ordinary notation. 
Sometimes the \dd Fimage $\ans{F(x):x\in X}$ of a set $X$ will be 
denoted by $F\ima X$. 

$\rV$ will denote the universe of all sets, 
$\rL$ the constructible universe. 

We shall use sans serif characters like $\rG$ to denote graphs 
and other binary relations. On the other hand, the math italic 
$G$ will denote generic sets. 

By {\it dense\/} we shall always mean: {\it open dense\/}.

\np

\subsection{Solovay and collapse models}
\label{somo}

This section reviews notation related to and some properties 
of the Solovay model and collapse generic models. (The latter 
type will be used in the proof of Theorem~\ref{tlast}.) 

For an ordinal $\al,$ $\smf\al=\col\al=\bigcup_{n\in\om}\al^n$ 
denotes the forcing to collapse $\al$ down to $\om.$ It is 
ordered as follows: $u\<v$ ($u$ is {\it stronger\/} than $v$) 
iff $v\sq u$.

For $\la\in\Ord,$ we let $\smfl\la$ ($\smfu\la$ in 
Solovay~\cite{sol}) be the ``finite support'' product of all 
sets $\smf\al,\msur$ $\al<\la.$ In other words, $\smfl\la$ is the 
set of all functions $p$ such that $\dom p$ is a finite subset  
of $\la$ and $p(\al)\in\smf\al$ for each $\al<\la,\msur$ 
$\al\in\dom p,$ with the order: $p\<q$ iff $\dom q\sq\dom p$ 
and $p(\al)\<q(\al)$ for all $\al\in\dom q$.  



\subsubsection{Solovay's ``\protect\dd{\protect\Sigma}lemma''}
\label{Sig}

It is known that the Solovay model is seen from any subclass of the 
form $\rL[x],$ where $x$ is a real or a set definable from a real, 
in one and the same regular way. This remarkable uniformity is 
partially inherited by collapse extensions. 
The following lemma stands behind this crucial property.  

\ble
\label{44}
Suppose that\/ $\gM$ is a transitive model of\/ $\ZFC,$ $P\in\gM$ 
is a \poset, and\/ $G\sq P$ is\/ \dd Pgeneric over\/ $\gM.$  Let\/ 
$S\in \gM[G]\,,\msur$ $S\sq\Ord.$ There exists a set 
$\Sg\sq P,\msur$ $\Sg\in \gM[S]$ such that\/ $G\sq\Sg$ and\/ 
$G$ is\/ \dd\Sg generic over\/ $\gM[S]$.
\ele
Thus $\gM[G]$ is a generic (not necessarily \dd Pgeneric) 
extension of any submodel $\gM[S],\msur$ $S\sq\Ord$.)\hspep

\proof{}(extracted from the proof of Lemma 4.4 in 
\cite{sol}). 

Let $\doS$ be a name for $S$ in the language of the 
\dd Pforcing  over $\gM$.

Define a sequence of sets 
$A_\al\sq P\;\;(\al\in\Ord)$ by induction on $\al$ in $\gM[S]$.
\ben
\def\theenumi{\hskip1pt(A\arabic{enumi})\hskip1pt}
\def\labelenumi{\theenumi}
\itemsep=1mm
\itla{aa1}\msur
$p\in A_0$ iff, for some $\sg\in\Ord,$ either $\sg\in S$ but $p$ 
\dd Pforces $\sg\not\in\doS$ over $\gM,$ or 
$\sg\not\in S$ but $p$ \dd Pforces $\sg\in\doS$ over $\gM$.

\itla{aa2}\msur
$p\in A_{\al+1}$ iff there exists a dense in $P$ below $p$ set 
$D\in \gM,\msur$ $D\sq A_\al$.

\itla{aa3}
If $\al$ is a limit ordinal then $A_\al=\bigcup_{\ba<\al}A_\ba$.\its
\een

The following properties of these sets are easily verifiable: 
first, if $p\in A_\al$ and $q\< p$ in 
$P$ then $q\in A_\al\,;$ second, if $\ba<\al$ then $A_\ba\sq A_\al$. 

Since each $A_\al$ is a subset of $P,$ it follows that 
$A_\da= A_{\da+1}$ for an ordinal $\da\in\gM.$ We put 
$\Sg=P\setminus A_\da.$ Thus $\Sg$ can be seen as the set of 
all conditions $p\in P$ which do not force something about $\doS$ 
which contradicts a factual information about $S.$ 
We prove that $\Sg$ is as required. This 
involves two auxiliary facts.
\ben
\def\theenumi{$(\Sigma\arabic{enumi})$}
\def\labelenumi{\theenumi\msur}
\itemsep=1mm
\itla{sig1}
$G\sq\Sg$. 

\itla{sig2}\psur
{\it If\/ $D\in \gM$ is a dense subset of\/ $P$ then\/ 
$D\cap\Sg$ is dense in\/ $\Sg$}.
\een
To prove \ref{sig1} assume on the contrary that 
$G\cap A_\ga\not=\emps$ for some $\ga.$ Let $\ga$ be the least 
such an ordinal. Clearly $\ga$ is not limit and $\ga\not=0;$ 
let $\ga=\al+1.$ Let ${p\in A_\ga\cap G.}$ Since $G$ is generic, 
item~\ref{aa2} implies ${G\cap A_\al\not=\emps},$ 
contradiction. \ref{sig2} is easy: if $p\in\Sg$ then 
$p\not\in A_{\da+1};$ therefore by~\ref{aa2} $D\cap A_\da$ 
is {\it not\/} dense in $P$ below $p,$ so there exists 
$q\<p,\msur$$q\in D\setminus A_\da$ -- then $q\in\Sg$. 

We prove that $G$ is \dd\Sg generic over $\gM[S].$ 
Suppose, towards contradiction, that $D\cap G=\emps$ for 
a dense in $\Sg$ set $D\in\gM[S],\msur$ $D\sq\Sg$. 

Since ${D\in\gM[S]},$ there exists an \dd\in formula $\Phi(x,y)$ 
containing only sets in $\gM$ as parameters and such that $\Phi(S,y)$ 
holds in $\gM[S]$ iff $y=D$.

Let $\Psi(G')$ be the conjunction of the following formulas:
\ben
\itemsep=1mm
\def\theenumi{\hskip1pt(\arabic{enumi})\hskip1pt}
\def\labelenumi{\theenumi}
\itla{1)}
\msur $S'=\doS[G']$ (the interpretation of the ``term'' $\doS$ via 
$G'$)\hspace{1mm}{} is a set of ordinals, and there exists unique 
$D'\in\gM[S']$~\footnote
{\ Let us assume, to avoid irrelevant complications, that $\gM$ 
satisfies the following axiom: $\exists\,X\,(\rV=\rL[X])\,;$ this 
makes $D'\in\gM[S']$ a legitimate formula. The lemma will be applied 
below only in this particular case, but the proof could be modified 
to fix the general case.} 
such that $\Phi(S',D')$ holds in $\gM[S']$;

\itla{2)}
this set $D'$ is a dense subset of $\Sg'$ where $\Sg'\sq P$ is the 
set obtained by applying our definition of $\Sg$ for $S=S'$ within 
$\gM[S']$;

\itla{3)}\msur
$D'\cap G'=\emps$.
\een 
Then $\Psi(G)$ is true in $\gM[G]$ by our assumptions. Let 
$p\in G$ \dd Pforce $\Psi(G)$ over $\gM.$ Then $p\in\Sg$ by 
\ref{sig1}. By the density there exists a condition 
$q\in D\,,\msur$ $q\<p\,.$ Les us consider 
a \dd\Sg generic over $\gM[S]$ set $G'\sq\Sg$ containing $q.$ 
Then $G'$ is also \dd Pgeneric over $\gM$ by \ref{sig2}.   
We observe that $\doS[G']=S$ because $G'\sq\Sg.$ It follows that 
$D'$ and $\Sg'$ (as is the description of $\Psi$) coinside with 
resp.\ $D$ and $\Sg.$ In particular ${q\in D'\cap G'},$ a 
contradiction because $p$ forces \ref{3)}.\qed

\bcor
\label{444+}
If\/ $P=\smf\la$ in the lemma for an ordinal\/ $\la\in\gM$ then\/ 
$\gM[G]$ is\/ \dd{\smf\al}generic extension of\/ $\gM[S]$ for an 
ordinal\/ $\al\<\la$.
\ecor
It is allowed that $\al=1,$ to include in particular the case 
$\gM[S]=\gM[G].$ Note that $\al,$ generally speaking, depends on 
$G:$ different parts of $\Sg$ may have different branching 
poperties.\hspep

\proof We observe that $\Sg$ is, in $\gM[S],$ a subtree of the full 
\dd\la branching tree $\smf\la=\col\la.$ This allows to derive the 
result by ordinary reasoning.\qed

\subsubsection{``Weak'' sets in the models}
\label{weak}

As usual, $\rL$ denotes the constructible universe. By 
\dd\Om{\em Solovay model axiom\/} and 
\dd\Om{\em collapse model axiom\/} 
we shall understand the hypotheses:
\bde
\itemsep=1mm
\item[$\sma\Om:$]
$\msur\Om$ is inaccessible in $\rL$ and the universe  
is a \dd{\smfl\Om}generic extension~of~$\rL.$ 

\item[$\clps\Om:$] 
$\msur\Om$ is a limit cardinal in $\rL$ and 
the universe is a \dd{\smf\Om}generic extension~of~$\rL.$ 
\ede
%
%
A set $x$ will be called \dd\Om{\em ``weak'' over\/ ${\gM}$\/} 
iff $x$ belongs to a \dd{\smf\al}generic extension of $\gM$ for 
some $\al<\Om.$ We define 
\dm 
\cW[S]=\ans{x\in\cD:x\,\hbox{ is \dd\Om ``weak'' over }\,\rL[S]},
\dm 
the set of all reals ``weak'' over a class $\rL[S].$ As 
usual $\cW=\cW[\emps]$.

\ble
\label{allw}
Assume\/ $\sma\Om.$ Then all reals are\/ \dd\Om``weak'' over\/ 
$\rL$.\qed
\ele
%
Lemma is not true in the assumption 
$\clps\Om,$ of course. This observation outlines the difference in 
the development of the Solovay and collapse models: usually properties 
of {\it all\/} reals in the former remain true for ``weak'' reals in 
the latter. 

\bpro
\label{1sm}
Assume\/ $\sma\Om$ {\rmt(}resp.\ assume\/ $\clps\Om${\rmt)}. Then\/ 
$\omi=\Om$  
{\rmt(}resp.\ $\omi=\Om^+,$ the next cardinal in\/ $\rL${\rmt)}.  
Furthermore, suppose that a set\/ $S\sq\Ord$ is\/ \dd\Om``weak''  
over\/ $\rL.$ Then\/ 
\ben
\def\itemsep{1mm}
\def\theenumi{{\arabic{enumi}}}
\def\labelenumi{{\rmt\theenumi.}}
\itla{sm1}\msur
$\Om$ is inaccessible in\/ $\rL[S]$ and the universe\/ 
$\rV$ is a\/ \dd{\smfl\Om}generic extension of $\rL[S].$ 
{\rmt(}Resp.\ $\Om$ is a limit  
cardinal in\/ $\rL[S]$ and\/ $\rV$ is a\/ \dd{\smf\Om}generic 
extension of\/ $\rL[S]$.{\rmt)}

\itla{sm2} 
If\/ $\Phi$ is a sentence containing only sets in\/ $\rL[S]$ as 
parameters then\/ $\La$ {\rm(}the empty function\/{\rm)} decides\/ 
$\Phi$ in the sense of\/ $\smfl\Om$ 
{\rmt(}resp.\ in the sense of\/ $\smf\Om${\rmt)} 
as a forcing notion over $\rL[S]$.

\itla{sm4} 
If a set\/ $X\sq\rL[S]$ is\/ $\od[S]$ then\/ $X\in\rL[S]$.
\een 
\epro
($\od[S]=S$\msur{}--{\it ordinal definable\/}, that is, 
definable by an \dd\in formula containing $S$ and ordinals as 
parameters.)\hspep

\proof {\em Item \ref{sm1}\/}. By definition, $S$ belongs to a   
\dd{\smf\al}generic extension of $\rL$ for some $\al<\Om.$ In the 
$\clps\Om$ case, the universe $\rV$ is a \dd{\smf\ga}generic 
extension of $\rL[S]$ for some $\ga\<\Om,$ by Corollary~\ref{444+}. 
But $\ga<\Om$ is impossible: indeed otherwise $\Om$ simply remains 
uncountable in $\rV,$ contradiction with $\clps\Om$.

Let us consider the $\sma\Om$ case. We have $\rV=\rL[G]\,;\msur$ 
$G$ being a \dd{\smfl\Om}generic set over~$\rL.$ Obviously 
$S\in\rL[x]$ for a real $x.$ It follows 
(Corollary 3.4.1 in \cite{sol}) that there exists an 
ordinal $\la<\Om$ such that $S$ belongs to $\rL[G_{\<\la}]$ 
for some $\la<\Om,$ where $G_{\<\la}=G\cap\cC{\<\la}$ is 
\dd{\smfe\la}generic over $\rL.$ Therefore by Lemma 4.3 in 
\cite{sol}, $\rL[G_{\<\la}]$ also is a \dd{\smf\la}generic 
extension of $\rL.$ Finally $\rL[G_{\<\la}]$ is a 
\dd{\smf\ga}generic extension of $\rL[S]$ for some 
$\ga\<\la$ by Corollary~\ref{444+}. 

On the other hand, $\rV=\rL[G]$ is a \dd{\smf{\geq\la+1}}generic 
extension of $\rL[G_{\<\la}],$ so that totally $\rV$ is a 
\dd{(\smf\ga\ti\smf{\geq\la+1})}generic extension of $\rL[S].$ But 
the product can be easily proved to include a dense subset order 
isomorphic to $\smfl\Om.$ This ends the proof of item 
\ref{sm1} of the proposition.

{\em Items \ref{sm2} and \ref{sm4}\/}. Easily follow from item 
\ref{sm1}.\qed

\np

\subsection{Reals in Solovay and collapse models}
\label{reals}

In this section, a useful coding of reals and sets of reals in 
generic models of the type we consider, is introduced. After this 
we study a topology on reals in the models, close to the 
Gandy -- Harrington topology.

\subsubsection{The coding}
\label{coding}

If a set $G\sq\smf\al$ is \dd{\smf\al}generic over a transitive 
model $\gM$ ($\gM$ is a set or a class) then $f=\bigcup G$ maps 
$\om$ onto $\al,$ so that $\al$ is countable in $\gM[G]=\gM[f].$ 
Functions $f:\om\,\lra\,\al$ obtained 
this way will be called \dd{\smf\al}{\em generic over\/ $\gM$.}

We let $\gfu\al(\gM)$ be the set of all \dd{\smf\al}generic over 
$\gM$ functions $f\in \al^\om.$ We further define 
$\gfu\al[S]=\gfu\al(\rL[S])$ 
and $\gfu\al=\gfu\al(\rL)=\gfu\al[\emps]$.

Recall that {\it reals\/} are elements of $\cD=2^\om$ in 
this research. 

Let $\al\in\Ord$. $\trb\al$ will denote the set of all ``terms'' --- 
indexed sets of the form $\zt=\ang{\al,\ang{\zt_n:n\in\om}}$ 
such that ${\zt_n\sq\smf\al}$ for each $n.$  
``Terms'' $\zt\in\trb\al$ can be used to code functions 
$C:\al^\om\;\lra\;\cD=2^\om;$ namely, for every $f\in\al^\om$ 
we define $x=\zC\zt(f)\in\cD$ by: $x(n)=1$ iff $f\res m\in \zt_n$ 
for some $m$.

We put $\trb{}=\bigcup_{\al<\Om}\trb\al.$ (This definition makes 
sense only when $\Om$ is fixed by the context, \eg\ provided 
we suppose $\sma\Om$ or $\clps\Om.$ Recall that $\Om=\omi$ 
assuming $\sma\Om$ but $\Om<\omi$ assuming $\clps\Om$.)  

Suppose that $\zt=\ang{\al,\ang{\zt_n:n\in\om}}\in\trb\al,\msur$  
$u\in\col\al=\smf\al,\msur$ $\gM$ arbitrary. We define  
$\zX{\zt u}(\gM)=\ans{\zC\zt(f):u\subset f\in\gfu\al(\gM)}$ and 
$\zX\zt(\gM)=\zX{\zt \La}(\gM)=\zC\zt\ima\gfu\al(\gM).$ In 
particular, let $\zX\zt[S]=\zX\zt(\rL[S])$ and 
$\zX\zt=\zX\zt[\emps]=\zX\zt(\rL);$ the same for $\zX{\zt u}$.

We recall that $\cW[S]$ is the set of all reals \dd\Om``weak'' 
over $\rL[S]$.

\bpro
\label{2sm}
Assume that either\/ $\sma\Om$ or\/ $\clps\Om$ holds in the 
universe. Let $S\sq\Ord$ be\/ \dd\Om``weak'' over\/ $\rL.$ Then
\ben
\itemsep=1mm
\def\theenumi{{\arabic{enumi}}}
\def\labelenumi{{\rmt\theenumi.}}
\itla{sm6} 
If\/ $\al<\Om,\msur$ $F\sq\gfu\al[S]$ is\/ $\od[S],$ and\/ 
$f\in F,$ then there exists\/ $m\in\om$ such that each\/ 
$f'\in\gfu\al[S]$ 
satisfying\/ $f'\res m= f\res m$ belongs to\/ $F$.

\itla{sm5} 
For each\/ $x\in\cW[S]$ there exist\/ $\al<\Om,\msur$ 
$f\in\gfu\al[S],$ and\/ $\zt\in\trb\al\cap\rL[S]$ such that\/ 
$x=\zC\zt(f)$.

\itla{xl1} 
Every\/ $\od[S]$ set $X\sq\cW[S]$ is a union of sets of the 
form\/ $\zX\zt[S],$ where $\zt\in\trb{}\cap\rL[S].$ Conversely\/ 
$\zX\zt[S]\sq\cW[S]$ whenever\/ $\zt\in\trb{}\cap\rL[S]$.

\itla{xl2} 
Assume that\/ ${\al<\Om},\msur$ ${u\in\col\al},$ and\/ 
${\zt\in\trb\al\cap\rL[S]}\,.$ Then each\/ $\od[S]$ set\/ 
$X\sq {\zX{\zt u}[S]}$ is a union of sets of 
the form\/ $\zX{\zt v}[S],$ where\/ $u\sq v\in\col\al.$
\een
\epro
(Take notice that $\cW[S]=\cD,$ all reals, provided $\sma\Om,$ 
but $\cW[S]$ is a proper subset of $\cD$ in the assumption of 
$\clps\Om$.)\hspep

\proof {\em Item \ref{sm6}\/}. We observe that  
$F=\ans{f'\in\al^\om:\Phi(S,f')}$ for an \dd\in formula $\Phi.$ 
Let $\Psi(S,f')$ denote the formula: \ 
``$\La$ \dd Pforces $\Phi(S,f')$ over the universe'', \ 
where $P$ is the forcing 
$\smfl\Om$ in the case $\sma\Om,$ and $\smf\Om$ otherwise, so  
\dm
F=\ans{f'\in\al^\om:\Psi(S,f')\,\hbox{ is true in }\,\rL[S,f']}
\dm
by Proposition~\ref{1sm} (items \ref{sm1}, \ref{sm2}). 
Thus, since $f\in F\sq\gfu\al[S],$ there exists $m\in\om$ such 
that the restriction $u=f\res m$\hspace{1mm}{} 
(then $u\in\col\al=\smf\al$)\hspace{1mm}{} 
\dd{\smf\al}forces $\Psi(S,\namef)$ over $\rL[S],$ 
where $\namef$ is the name of the \dd\al collapsing function. 
The $m$ is as required. 

{\em Item \ref{sm5}\/}. By definition $x$ belongs to an 
\dd{\smf\al}generic extension of $\rL[S],$ for some ordinal 
$\al<\Om.$ Thus $x\in\rL[S,f]$ where $f\in\gfu\al[S].$ Let 
$\namex$ be a name for $x.$ It remains to define 
$\zt_n=\ans{u\in\smf\al:u\,\hbox{ \dd{\smf\al}forces }\,
\namex(n)=1}$ in $\rL[S]$.

{\em Item \ref{xl1}\/}. Let $x\in X.$ We use item 2 to get an 
ordinal $\al<\Om,$ a generic function $f\in\gfu\al[S],$ and a 
``term'' $\zt\in\trb\al\cap\rL[S]$ such that\/ $x=\zC\zt(f).$ 
Applying item 1 to the set 
$F=\ans{f'\in\gfu\al[S]:\zC\zt(f')\in X}$  
and the given function $f,$ we obtain a forcing condition 
$u=f\res m\in\smf\al\msur$ $(m\in\om)$  
such that $x\in \zX{\zt u}[S]\sq X.$ Finally the set 
$\zX{\zt u}[S]$ is equal to $\zX{\zt'}[S]$ for some other 
$\zt'\in \trb\al\cap\rL[S]$.

{\em Item \ref{xl2}\/}. Similar to the previous item.\qed

\subsubsection{The ordinal definable topology}
\label{tt}
 
It occurs that the topology generated by $\od$ sets in the 
Solovay and collapse models has a strong semblance of the 
Gandy -- Harrington topology 
(in a simplified form because the specific $\is11$ technique  
becomes obsolete). In particular, the topology is strongly 
Choquet. However we shall not utilize this property. The 
treatment will be organized in a rather forcing--like way, 
in the manner of Miller~\cite{mill}.

 A set $X$ will be called \dd\od{\it 1st-countable\/} if the $\od$ 
power set $\oP(X)=\cP(X)\cap\od$ is at most countable. In this case, 
$\oP(X)$ has only countably many different $\od$ subsets 
because it is a general property of the Solovay model that 
$\oP(\cX)$ is countable for any countable $\od$ set 
$\cX\sq\od$.

\ble
\label{dizl}
Assume that either\/ $\sma\Om$ or\/ $\clps\Om$ holds in the 
universe. Let\/ $\al<\Om$ and\/ $\zt\in\trb\al\cap\rL.$ 
Then\/ $X=\zX\zt$ is\/ \dd\od 1st-countable.
\ele
\proof By Proposition~\ref{2sm} every $\od$ subset of $X$ is 
uniquely determined by an $\od$ subset of $\smf\al=\col\al.$ 
Since each $\od$ set $S\sq\smf\al$ is constructible 
(Proposition~\ref{1sm}), 
we obtain an $\od$ map $h:\al^+\,\hbox{ onto }\,\oP(X),$ 
where $\al^+$ is the least \dd\rL cardinal bigger than $\al.$ 
Therefore $\oP(X)$ has at most \dd{\al^{++}}many $\od$ subsets. 
It remains to notice that $\al^{++}<\Om$ follows from either of  
$\sma\Om$ or $\clps\Om$.\hspep\qed

Let  
$\dX=\ans{X\sq\cW:X\,\hbox{ is }\,\od\;\hbox{ and nonempty}\,},$ 
where $\cW=\cW[\emps]$. 

(We recall that $\cW=\cD=\hbox{all reals}$ 
in the assumption $\sma\Om$ but $\cW\sneq\cD$  
in the assumption $\clps\Om,$ consisting of all reals 
\dd\Om``weak'' over $\rL$.) 

Let us consider $\dX$ as a forcing notion (smaller sets are 
stronger conditions). We say that a set $G\sq\dX$ is 
\dd\od{\it generic in} $\dX$ iff it nonempty intersects every 
dense $\od$ subset of $\dX$.~\footnote
{\ This can be transformed to an ordinary forcing over $\rL$ in the 
assumption $\sma\Om.$ Of course formally $\dX\not\in\rL,$ but 
$\dX$ is $\od$ order isomorphic to a partially ordered set in $\rL$.  
}

\ble 
\label{sing}
Assume that\/ $\sma\Om$ or\/ $\clps\Om$ holds. If a set\/ $G\sq\dX$ 
is\/ \dd\od gen\-eric then the intersection\/ $\bigcap G$ is a 
singleton $\ans{x}=\ans{x_G}$ and\/ $x_G\in\cW$.
\ele
\proof It follows from the density and Proposition~\ref{2sm} 
that $G$ contains a set $X$ of the form $X=\zX\zt$ for some 
$\zt\in\trb\al,\;\,\al<\Om.$  Note that by Proposition~\ref{2sm} 
(item~\ref{xl2}) for any dense in $\smf\al=\col\al$ $\od$ set 
$D\sq\smf\al$ and any $u\in\col\al$ the set 
\dm
\cX=\ans{Y\in\dX:X\cap Y=\emps}\cup 
\ans{\zX{\zt v}: v\in D\cj u\sq v} 
\dm
is $\od$ and dense in $\dX.$ This leads to a function 
$f\in\gfu\al$ such that $\zX{t\;f\res m}\in G$ for all $m\in\om,$
since the family of all $\od$ subsets of $\smf\al$ is countable 
(see above).  

Let us show that $x=\zC\zt(f)\in\bigcap G.$ Let $Y\in G.$ 
The set $D$ of all $u\in\smf\al$ such that either 
$\zX{\zt u}\sq Y$ or $\zX{\zt u}\cap Y=\emps$ is dense in 
$\smf\al$ by Proposition~\ref{2sm} (item~\ref{xl2}), and 
obviously $\od.$ Therefore there exists a 
number $m$ such that $u=f\res m\in D.$ Then $\zX{\zt u}\in G$  
by the choice of $f,$ so that $\zX{\zt u}\cap Y=\emps$ is 
impossible. Thus $\zX{\zt u}\sq Y$ and $x=\zC\zt(f)\in Y,$ as 
required.\qed

\np

\subsection{The dichotomy}
\label{1st}

Let us start the proof of the dichotomy of Theorem~\ref{main}
in the assumption $\sma\Om.$ In parallel, consideration of 
the collapse (\ie\ \dd{\smf\Om}generic) 
models will be carried out, in order to prepare the 
proof of the dichotomy of Theorem~\ref{tlast} afterwards.  

This line of arguments takes this section (mainly devoted to the 
case when a colouring exists) and two next sections (where we 
define a homeomorphism or embedding of $\Go$ into the gived graph). 
During the reasoning we shall assume that the following hypotheses 
is satisfied:
\ben
\itemsep=1mm
\def\theenumi{$(\ddag)$}
\def\labelenumi{\theenumi}
\itla{h}
either $\sma\Om$ or $\clps\Om$ holds, and $\rG$ is an $\od$ graph 
on reals satisfying $\rG\sq\cWt,$ where $\cWt$ is the 
two--dimentional copy of $\cW$, 
\dm
\cWt=\ans{\ang{x,y}\in\cD^2:
\ang{x,y}\,\hbox{ is \dd\Om``weak'' over }\,\rL}\,.
\dm
\een 
The $\sma\Om$ part of \ref h will work towards the proof of 
Theorem~\ref{main}; note that $\cWt=\cD\ti\cD$ in the assumption 
$\sma\Om$ (all reals and pairs of reals are ``weak'').  

The $\clps\Om$ part will be involved in the proof of 
Theorem~\ref{tlast}. Notice that $\cWt\sq\cW\ti\cW$ but it is not 
clear whether $\cWt\sneq\cW\ti\cW$ in the assumption $\clps\Om.$ 
(We do not know how to find a non--``weak'' pair of ``weak'' reals.) 

\subsubsection{Incompatibility}
\label{incomp}

We prove first of all that the requirements \ref{1} and \ref{2} 
of Theorem~\ref{main} are incompatible in the assumption of $\sma\Om.$  
Suppose on the contrary that $\rG$ satisfies both \ref{1} and \ref{2}. 
Then $\Go$ itself admits a $\rod$ colouring $c:\cD\,\lra\,\omi$. 

We get the contradiction following an argument for Borel 
colourings in \cite{kst}. Each set $X_\al=c^{-1}(\al)$ ($\al<\omi$) 
is an $\od$ set of reals, 
and $\cD=\bigcup_{\al<\omi} X_\al.$ By the known properties 
of the Solovay model, some $X_\al$ is not meager. Then $X_\al$ is 
co-meager on a set of the form $\cD_w=\ans{x\in\cD:w\sq x},$ where 
$w\in 2^\lom.$ By the choice of the sequence of tuples 
$\ww_n\in 2^\lom$ (see Introduction), we have $w\sq\ww_n$ for some 
$n.$ Then $X_\al$ is co-meager on  $\cD_{\ww_n}.$ Let $H$ be the 
automorphism of $\cD$ defined by: $H(x)=y,$ where $y(k)=x(k)$ for all 
$k\not=n,$ but $y(n)=1-x(n).$ Then $Y=H^{-1}(X_\al\cap \cD_{\ww_n})$ is 
co-meager on $\cD_{\ww_n}.$ 
Let $x\in X_\al\cap \cD_{\ww_n}\cap Y.$ Then $y=H(x)\in X_\al,$ 
but $x\Go y,$ contradiction. 

\subsubsection{The dichotomy}

We recall that a set $X$ is \dd\rG{\it discrete\/} iff 
$G\cap (X\ti X)=\emps.$ 
The following theorem will imply the dichotomy results of 
theorems \ref{main}, \ref{spe}, and \ref{tlast}.

\bte
\label{m2}
\label{m1} 
Assume\/ \ref h. Then 
\ben
\itemsep=1mm
\def\theenumi{\hskip1pt{(\roman{enumi})}\hskip1pt}
\def\labelenumi{{\rmt\theenumi}}
\itla{mm1}\msur 
If \underline{every} real\/ $x\in\cW$ belongs to an\/ $\od$ 
\dd\rG discrete set then\/ $\rGw$ admits an\/ $\od$ colouring by 
countable ordinals. 

\itla{mm2}
If \underline{not every} real\/ $x\in\cW$ belongs to an\/ 
$\od$ \dd\rG discrete set then\/ $\Go\hoc \rGw.$ Moreover, we 
have\/ ${\Go\emc\rGw}$ in each of the three cases$:$ 
\ben
\itemsep=1mm
\def\labelenumii{{\rmt(\theenumii)}}
\itla{2a} $\rGw$ is locally countable$;$ 

\itla{2b} $\rGw$ is acyclic$;$

\itla{2c} $x\rGw y$ implies\/ $x\not\in\rL[y]\cj y\not\in\rL[x]$.  
\een
\een
\ete
We shall demonstrate, after quite a short proof of statement 
\ref{mm1}, how this theorem implies theorems \ref{main} and 
\ref{spe}. Statement~\ref{mm2} of the theorem will need much more 
efforts; the proof takes the next two sections.\hspep
  
\proof{}of item \ref{mm1} of Theorem~\ref{m1}. In principle the 
assumption immedialely leads to an $\od$ colouring of $\rG$ by 
ordinals. Indeed there exists an \dd\in definable map 
$\phi:\Ord\,\hbox{ onto }\,\od.$ For each real $x,$ we let $c(x)$ be 
the least ordinal $\ga$ such that $x\in\phi(\ga)$ and $\phi(\ga)$ is 
\dd\rG discrete. Then $c$ is an $\od$ colouring of $\rG$ by ordinals.

A more difficult problem is to get a colouring by 
{\it countable\/} ordinals. 

The set $\trb{}\cap\rL$ is obviously a constructible set of 
cardinality $\Om$ in $\rL.$ 
(Notice that $\Om=\omi$ in the case $\sma\Om$ and $\Om<\omi$ in 
the case $\clps\Om$.) Let us fix an \dd\in definable 
enumeration $\trb{}\cap\rL=\ans{\zt_\al:\al<\Om}.$ It follows 
from Proposition~\ref{2sm} that every $\od$ set of reals is 
a union of sets of the form $Z_\al=\zX{\zt_\al}.$ Now one 
defines an $\od$ colouring $c$ of $\rG$ by countable ordinals, 
setting $c(x)$ to be the least $\al<\omi$ such that 
$x\in Z_\al$ and $Z_\al$ is a \dd\rG discrete set.\hspep\qedd

\noi{\bf Remark.} 
Let $\rG$ be a $\is1n$ graph $(n\geq2)$ in Theorem~\ref{m2}, 
item~\ref{mm1}, in the $\sma\Om$ version. Then $\rG$ admits a\/ 
$\idh n$ colouring. Indeed one easily sees that the relation \ 
``$\zX\zt$ is \dd\rG discrete'' \ is $\ip1n,$ hence 
$\iph{n-1}.$ The enumeration of terms can be chosen to be 
$\idh1.$ An elementary computation shows that then  
$c$ is $\idh n$.

\subsubsection{Applications}
\label{apli}

\proof{}of Theorem~\ref{main}. Consider a $\rod$ graph $\rG$ on 
reals in the assumption $\sma\Om.$ We shall assume that $\rG$ 
is in fact an $\od$ graphs. (The ``boldface'' case, when $\rG$ is 
$\od[p]$ for a real $p,$ can be handled similarly: $p$ uniformly 
enters the reasoning as a parameter. In particular $\rL$ changes 
to $\rL[p]$ in the definition of $\cW,$ similarly $\od$ changes 
to $\od[p]$ in the definition of the topology and the set $\dX$ 
in Subsection~\ref{tt}, \etc) 

Since clearly $\cW=\cD,$ all reals, in the assumption $\sma\Om,$ 
we have \ref h. This allows to apply Theorem~\ref{m2}, immediately 
getting the dichotomy of Theorem~\ref{main} for the given 
graph $\rG$.\hspep\qed

\proof{}of Theorem~\ref{spe}. Let us consider, in the assumption 
$\sma\Om,$ a $\rod$ graph $\rG.$ As above, we shall suppose 
that $\rG$ is actually an $\od$ graph. Let us divide $\rG$ 
into the two parts: 
${\rG'}=\ans{\ang{x,y}\in{\rG}:x\in \rL[y]
\hspace{1pt}\hbox{ or }\hspace{1pt}y\in \rL[x]},$
and ${\rG''}={\rG}\setminus{\rG'}$. 

It is asserted that $\rG'$ is a special $\od$ graph.

Indeed we have to present $\rG'$ as the union of an $\od$ 
\dd{\omi}sequence of curves. Let, for $\al<\omi\,,$ 
${\rK}_{\al0}$ be the set of all pairs $\ang{x,y}\in{\rG}$ 
such that $x$ is the \dd\al th real in the sense of the 
\dd yconstructible wellordering of $\rL[y].$ Let 
${\rK}_{\al1}$ be the set of all pairs $\ang{x,y}\in{\rG}$ 
such that $y$ is the \dd\al th real in the sense of the 
\dd xconstructible wellordering of $\rL[x].$ Then by 
definition $\rG'$ is the union of all these sets, all of 
them are curves, and they can be ordered in an $\od$ 
\dd{\omi}sequence.

Thus by the assumption of Theorem~\ref{spe} the special 
subgraph $\rG'$ satisfies the dichotomy \ref1 \vs\ 
${\Go\emc\rG}$ of Theorem~\ref{main}. 

On the other hand $\rG'',$ the other subgraph, 
satisfies \ref{2c} of Theorem~\ref{m2} by definition. It follows 
that $\rG''$  satisfies the same dichotomy, as required. \qed

\np

\subsection{The homomorphism}
\label{2nd}

In this section, we present the proof of item~\ref{mm2} of 
Theorem~\ref{m2} in the version $\Go\hoc\rG.$ 
Thus it is assumed that either $\sma\Om$ or $\clps\Om$ holds. We 
consider an $\od$ graph $\rG$ on reals, and suppose that the set 
\dm
A_0=\ans{x\in\cW:x\,
\hbox{ does not belong to an $\od$ \dd\rG discrete set}}
\dm
is nonempty. A continuous homomorphism $H$ of $\Go$ into $\rG$ 
will be defined. 

Sinse $A_0$ is obviously an $\od$ set, there exists, by 
Proposition~\ref{2sm}, a set $A\sq A_0$ of the form 
$A=\zX\zt,$ where $\zt\in\trb\al\cap\rL$ for some $\al<\Om$. 

It follows from Lemma~\ref{dizl} that the set 
$\dX_A=\ans{X\in \dX:X\sq A}$ contains only countably many 
$\od$ subsets in the universe. Let $\ans{\cX_n:n\in\om}$ be a 
(not necessarily $\od$) enumeration of all dense subsets of 
$\dX_A$.

In general for any \dd\od 1st-countable set $Q$ let us fix 
once and for all an enumeration ${\ans{\cX_n(Q):n\in\om}}$ of 
all dense $\od$ subsets 
in $\oP(Q)\setminus\ans{\emps},$ where, we recall, 
${\oP(Q)=\cP(Q)\cap\od}$ is the $\od$ part of the power set 
of $Q.$ No uniformity on $Q$ is assumed in the choice of 
the enumerations.

We recall that a sequence of binary tuples $\ww_n\in 2^n$ is 
fixed by the definition of $\Go,$ see Introduction. 
Let $m\in\om.$ By a {\it crucial pair\/} in $2^m$ we shall 
understand any pair $\ang{u,v}$ of tuples $u,\,v\in 2^m$ 
such that $u=\ww_k\we 0\we w$ and $v=\ww_k\we 1\we w,$ 
for some $k<m$ and $w\in 2^{m-k-1}$.  

We observe that the (directed) graph of crucial pairs in $2^m$ is 
a {\it tree\/}: each pair of $u,\,v\in 2^m$ is connected in $2^m$ 
by a unique (non--self--intersecting) chain of crucial pairs 
$u=u_0\spp u_1\spp u_2\spp ...\spp u_{n-1}\spp u_n=v,$ where 
$u'\spp v'$ means that either $\ang{u',v'}$ or $\ang{v',u'}$ is 
a crucial pair.

We shall define a set $X_u\sq A$ for every $u\in 2^\lom,$ 
and a binary relation $\rQ{uv}$ for each crucial pair 
$\ang{u,v}$ satisfying the following requirements:
\ben
\def\theenumi{C-\arabic{enumi}}
\def\labelenumi{\theenumi:}
\itemsep=1mm
\itla{x1}
each $X_u\;\,(u\in 2^m)$ is an $\od$ nonempty subset of $A$ and 
$X_u\in\cX_m$; 

\itla{x2}\msur
$X_{u\we i}\sq X_u$ for all $u\in 2^\lom$ and $i=0,\,1$.

\itla{j3}
if $\ang{u,v}$ is a crucial pair then $\rQ{uv}\sq \rG$ and 
$\rQ{uv}$ is a nonempty $\od$ set$;$

\itla{j2}
if $\ang{u,v}$ is a crucial pair then 
$\rQ{u\we i\,,\,v\we i}\sq\rQ{uv}$ for any $i=0,\,1$; 

\itla{j4}
if $\ang{u,v}$ is a crucial pair then $X_u \rQ{uv} X_v,$ that 
is,\\[1ex]
$\forall\,x\in X_u\;\exists\,y\in X_v\;(x\rQ{uv} y)
\hspace{5mm}\hbox{and}\hspace{5mm}
\forall\,y\in X_v\;\exists\,x\in X_u\;(x\rQ{uv} y)$;

\itla{j0} if $m\geq1$ and $u=\ww_{m-1}\we 0,\msur$ $v=\ww_{m-1}\we 1$ 
(a special type of crucial pairs) 
then $\rQ{uv}$ is an \dd\od 1st-countable set$;$

\itla{j1}
for any crucial pair of $u=\ww_n\we 0\we w$ and 
$v=\ww_n\we 1\we w$ in $2^m\msur$ $(n<m\in\om),$ 
$\rQ{uv}\in \cX_m(\rQ{\ww_n\we 0,\,\ww_n\we 1})$.
\een
We can assume without any loss of generality that 
$\rQ{uv}\sq X_u\ti X_v$ for all crucial pairs (for if not 
replace the relations by the intersections with 
$X_u\ti X_v$).

Suppose that such a system of sets $X_u$ and relations $\rQ{uv}$ 
has been defined. It follows from \ref{x1} and \ref{x2} that 
every sequence of the form $\ang{X_{a\res m}:m\in\om}$ (where 
$a\in\cD=2^\om)$) is \dd\od generic in $\dX,$ therefore the 
intersection $\bigcap_m X_{a\res m}$ is a singleton, say $H(a),$ 
by Lemma~\ref{sing}, and $H$ continuously maps $\cD$ into $\cD.$ 

Let us prove that $H$ is a homomorphism $\Go$ into $\rGw,$ 
that is, ${a\Go b}$ implies ${H(a)\rGw H(b)}$ for all 
$a,\,b\in\cD.$ By definition ${a\Go b}$ means that, for some 
$k,$ we have for instance $a=\ww_k\we 0\we z$ and 
$b=\ww_k\we 1\we z$ where $z\in 2^\om.$ It follows from 
\ref{j0}, \ref{j1}, \ref{j2} that the sequence of sets 
$Q_m=\rQ{\ww_k\we 0\we(z\res m),\,\ww_k\we 1\we(z\res m)}$ 
$\msur (m\in\om)$ is \dd\od generic in the two--dimentional 
``copy'' of $\dX,$ so that the intersection 
$\bigcap_m Q_m$ is a singleton. It follows from \ref{j4} that 
the singleton is equal to $\ang{H(a),H(b)}.$ We conclude that 
$\ang{H(a),H(b)}\in\rGw$ because 
$Q_0=\rQ{\ww_k\we 0,\,\ww_k\we 1}\sq \rGw$ by \ref{j3}. Thus 
$H(a)\rGw H(b)$ as required.

\subsubsection*{Construction of the sets and relations}

To begin with, we note that, since $A$ (see above) is a 
nonempty $\od$ set, it contains a subset $X\sq A,\msur$ 
$X\in\cX_0.$ We take such a set $X$ as $X_\La$. 

Suppose that $m\geq1,$ and the construction of sets 
$X_s$ and relations $\rQ{st}$ 
(where $\ang{s,t}$ is a crucial pair) has been completed 
for all $s,\,t\in 2^{<m}.$ Let us define the sets $X_u$ and 
relations $\rQ{uv}$ for $u,\,v\in 2^m,$ the next level. 

At the beginning, we set $X_{s\we i}=X_s$ for all $s\in 2^{m-1}$ 
and $i=0,\,1,$ and $\rQ{s\we i\,,\,t\we i}=\rQ{st}$ for all 
$i=0,\,1,$ and any crucial pair $\ang{s,t}$ in $2^{m-1}.$ This 
definition gives the initial versions of sets $X_u\msur$ $(u\in 2^m)$ 
and relations $\rQ{uv}$ ($\ang{u,v}$ being a crucial pair in $2^m,$ 
still with the exception of the crucial pair of 
$u_0=\fsg_{m-1}\we 0$ and $u_1=\fsg_{m-1}\we 1$.) The sets 
$X_u$ and relations $\rQ{uv}$ will be reduced, in several steps, 
to meet requirements \ref{x1} through \ref{j1}. 

{\it Step 1\/}.
This step is devoted to the ``new'' crucial pair of  
$u_0=\fsg_{m-1}\we 0$ and $u_1=\fsg_{m-1}\we 1$ in $2^m.$ At the 
moment $\rQ{u_0u_1}$ has not yet been defined. 

The set $X=X_{u_0}=X_{u_1}=X_{\fsg_{m-1}}$ is obviously a nonempty 
$\od$ subset of $A$. 

Notice that the $\od$ set ${\rC}={\rGw}\cap (X\ti X)$ is nonempty   
because $X$ is not \dd\rG discrete (being a subset of $A_0$). Using 
Proposition~\ref{2sm} (item~\ref{xl1}) and Lemma~\ref{dizl} (in 
the two--dimentional version), we get a nonempty \dd\od 1st-countable 
$\od$ set ${\rQ{}}\sq{\rC}.$ We finally put 
${\rQ{u_0u_1}}={\rQ{}}$ and $B_0=\dom \rQ{},\msur$ $B_1=\ran\rQ{},$ 
so that $B_0$ and $B_1$ are $\od$ subsets of $X_{\fsg_{m-1}},$ 
$\rQ{u_0 u_1}$ is an $\od$ relation and an \dd\od 1st-countable set  
(so \ref{j0} is satisfied), and $B_0\rQ{u_0 u_1}B_1$.\vom

{\it Step 2\/}. We reduce the current system of sets $X_u$ to 
the system of still $\od$ and nonempty sets $B_u\sq X_u$ 
such that $B_{u_0}=B_0,\msur$ $B_{u_1}=B_1,$ and \ref{j4} is 
satisfied. Let $u\in 2^m$ be other than $u_0$ or $u_1.$ Then 
$u$ can be connected in $2^m$ with some $u_i$ ($i=0,\,1$) 
by a unique non-self-intersecting chain of 
crucial pairs  
$u_i=v_0\spp v_1\spp v_2\spp...\spp v_n=u$ so that $u_{1-i}$ 
does not appear in the chain. 
We define $B^0=B_i$ and then, by induction, 
$
{B^{l+1}=\ans{y\in X_{v_{l+1}}:\exists\,x\in B^l\;(x\rQ{l}y)}}
$ 
for all $l<n,$ where ${\rQ{l}}={\rQ{v_lv_{l+1}}}$ or resp.\ 
${\rQ{l}}={(\rQ{v_{l+1}v_l})^{-1}}$ provided $\ang{v_l,v_{l+1}}$ 
or resp.\ $\ang{v_{l+1},v_l}$ is a crucial pair. 
The final set $B^n$ can be taken as $B_u$.\vom

{\it Step 3\/}. 
We have a system of $\od$ nonempty sets $B_u\;\,(u\in 2^m)$ and 
$\od$ relations $\rQ{uv}$ ($\ang{u,v}$ being a crucial pair in $2^m$) 
which satisfies all of \ref{x1} through \ref{j0}, still with the 
exception of \ref{j1} and the requirement $X_u\in\cX_m$ in \ref{x1}.  

To satisfy the latter let us fix $u\in 2^m.$ By the density of 
$\cX_m$ there exists a nonempty $\od$ set $B\sq B_u$ such that 
$B\in\cX_m.$ Using the construction of Step 2, we obtain 
a system of $\od$ nonempty sets $B'_u\sq B_u$ keeping the 
mentioned properties of the sets $B_u$ and in addition 
satisfying $B'_u=B$. 
 
We iterate this reduction consecutively $2^m$ times (the number 
of elements in $2^m$), getting a system of nonempty $\od$ sets   
$C_u\sq B_u$ which satisfies all of \ref{x1} through \ref{j0}, 
still with the exception of \ref{j1}. Let us take intersections 
of the form ${{\rQ{uv}}\cap(X_u\ti X_v)}$ as the ``new'' 
$\rQ{uv}$.\vom

{\it Step 4\/}.  
At this step, further reductions are executed, to fulfill 
\ref{j1}. We consider a crucial pair of 
$u'=\ww_n\we 0\we w$ and $v'=\ww_n\we 1\we w\msur$ $(n<m)$
in $2^m.$ Let the ``new'' $\rQ{u'v'}$ be an $\od$ 
subset of the ``old'' $\rQ{u'v'}$ which belongs to 
$\cX_m(\rQ{\ww_n\we 0,\,\ww_n\we 1}).$ 
The ``new'' sets $C'_{u'}$ and $C'_{v'}$ are defined by 
\dm 
C'_{u'}=\ans{x\in C_{u'}:\exists\,y\in C_{v'}\,({x\rQ{u'v'} y})}
\,,\hspace{3mm}
C'_{v'}=\ans{y\in C_{v'}:\exists\,x\in C_{u'}\,({x\rQ{u'v'} y})}\,,
\dm
so that still $C'_{u'}\rQ{u'v'} C'_{v'}.$ Then we restrict the 
sets $C_u\;\,(u\in 2^m)$ using the method described at Step 2. 
Running this construction consecutively for all crucial 
pairs $\ang{u',v'}$ in $2^m,$ we obtain finally a system of sets 
$X_u$ and relations $\rQ{uv}$ which satisfies all of \ref{x1} 
through \ref{j1}.\vom

This ends the inductive construction of sets $X_u$ and 
relations $\rQ{uv},$ and completes the proof of statement \ref{mm2} 
of Theorem~\ref{m2} in the version $\Go\hoc \rG$.\qedd

\np

\subsection{The embedding}
\label{3rd}

This section proves the ``moreover'' part of item~\ref{mm2} of 
Theorem~\ref{m2} where ${\Go\emc\rG}$ is requested. We shall also 
prove Theorem~\ref{spe}. We still suppose \ref h,  
consider an $\od$ graph $\rG,$ and continue to assume that the set 
$A_0$ (see the beginning of Section~\ref{2nd}) is nonempty.

The common denominator is to modify the splitting construction 
described in Section~\ref{2nd} so that the resulting function 
$H$ is a $1-1$ map. The following additional requirement is 
clearly sufficient. 
\ben
\def\theenumi{C-\arabic{enumi}}
\def\labelenumi{\theenumi:}
\itemsep=1mm
\addtocounter{enumi}{7}
\itla{x3}
if $s\in 2^{m-1}\;\;(m\geq 1)$ then 
$X_{s\we 0}\cap X_{s\we 1}=\emps$. 
\een
Let us demonstrate how the construction of Section~\ref{2nd} 
can be improved to meet this additional requirement, in the 
subcases \ref{2a}, \ref{2b}, \ref{2c} of Theorem~\ref{m2}. The 
result for \ref{2c} will be used to accomplish the proof of 
Theorem~\ref{spe}.

\subsubsection{Acyclic case}
\label{acy}

Let us assume that $\rGw$ is an acyclic graph. The 
construction of Section~\ref{2nd} expands by one more 
step.\vspace{1mm}

{\it Step 5\/}. At the moment we have a system of sets $X_u$ 
and relations $\rQ{uv}$ which satisfies all of \ref{x1} through 
\ref{j1}. Let us re--denote the sets, putting $D_u=X_u.$ The 
aim is to reduce sets $D_u\;\,(u\in 2^m)$ once again, to  
satisfy also \ref{x3}. 


Suppose that $s\in 2^{m-1}$ and $u_0=s\we 0,\;u_1=s\we 1.$ 
Then $u_0$ can be connected 
with $u_1$ in $2^m$ by a unique sequence of crucial pairs 
$u_0=v_0\spp v_1\spp ...\spp v_n=u_1$ where, by the choice of 
$u_0$ and $u_1,$ $n$ is necessarily an odd number. It follows from 
\ref{j4} that there exists a system of points $x_k\in D_{v_k}$ 
such that $x_k\rQ{k}x_{k+1}$ for all $k<n,$ where 
$\rQ{k}=\rQ{v_kv_{k+1}}$ or resp.\ $\rQ{k}=(\rQ{v_{k+1}v_k})^{-1}$ 
provided $\ang{v_k,v_{k+1}}$ or resp.\ $\ang{v_{k+1},v_k}$ is a 
crucial pair. 

Take notice that $x_0\not=x_n$ because $\rGw$ is acyclic and $n$ 
is odd, in particular $n\not=2.$ Therefore there exist disjoint 
$\od$ sets $D'_i\sq D_{u_i}\,\;(i=0,\,1)$ such that $x_0\in D'_0$ 
and $x_n\in D'_1.$ Now the construction of Step 2 in 
Section~\ref{2nd} allows to obtain a system of nonempty $\od$ 
sets $X_u\sq D_u$ $\msur (u\in 2^m)$ which keeps \ref{x1} through 
\ref{j1} and satisfies $X_{u_0}\sq D'_0$ and $X_{u_1}\sq D'_1,$ 
so $X_{u_0}\cap X_{u_1}=\emps$. 

Iterating this procedure for all relevant pairs in $2^m,$ we 
get finally a system at level $m$ satisfying all of \ref{x1} 
through \ref{j1}, and in addition \ref{x3}, as required. 

\subsubsection{Locally countable case}
\label{lc}

The case when $\rGw$ is a locally countable graph needs 
another approach. We let a 
{\it {\psur}1$\msur-\msur$1{\psur} curve\/} be a binary 
relation $\rK$ which is formally (the graph of) a $1-1$ 
map.

\ble
\label{11c}
Assume\/ $\sma\Om$ or\/ $\clps\Om\,.$ Let\/ $\rG$ be a locally 
countable\/ $\od$ graph. Then\/ $\rG$ is a union of\/ 
\dd{\aleph_1}many\/ $\od$ {\psur}1$\msur-\msur$1{\psur} curves.
\ele
\proof For any real $x,$ the set $G(x)=\ans{y:x\rG y}$ of all 
\dd\rG neighbours of $x$ is a countable $\od[x]$ set of 
reals. It is known from the theory of the Solovay and collapse 
models that then $G(x)\sq \rL[x].$ So 
$x\rG y$ implies $x\in \rL[y]$ and $y\in \rL[x].$ Thus  
${\rG}=\bigcup_{\al,\,\ba<\omi}{\rK}_{\al\ba}$ where 
${\rK}_{\al\ba}$ is the set of all pairs $\ang{x,y}\in{\rG}$ 
such that $x$ is the \dd\al th real in the sense of the 
\dd xconstructible wellordering of $\rL[y]$ and $y$ is the 
\dd\ba th real in the sense of the \dd yconstructible 
wellordering of $\rL[x]$.\hspep\qed

The lemma allows to carry out the construction in 
Section~\ref{2nd} so that, in addition to requirements 
\ref{x1} through \ref{j1} the following is satisfied:
\ben
\def\theenumi{C-\arabic{enumi}}
\def\labelenumi{\theenumi:}
\itemsep=1mm
\addtocounter{enumi}{8}
\itla{j6}
if $u=\ww_{m-1}\we 0,\;\,v=\ww_{m-1}\we 1,$ 
then $X_u\cap X_v=\emps$ and $\rQ{uv}$ is a $1-1$ curve. 
\een
Let us prove that this implies \ref{x3}. 
%
%
We have to check that 
$X_{s\we 0}\cap X_{s\we 1}=\emps$ for all $s\in 2^{m-1}.$ 
We recall that $\ww_{m-1}\in 2^{m-1}$ too, therefore as above 
there exists a unique chain of crucial pairs 
$\ww_{m-1}=s_1\spp s_2\spp ... \spp s_n=s$ in $2^{m-1}$. 

We prove by induction on $k$ that 
$X_{s_k\we 0}\cap X_{s_k\we 1}=\emps.$

The result for $k=1$ immediately follows from \ref{j6}. 

Let us suppose that $X_{s_k\we 0}\cap X_{s_k\we 1}=\emps$ and 
prove $X_{s_{k+1}\we 0}\cap X_{s_{k+1}\we 1}=\emps$.

One of the pairs $\ang{s_k,s_{k+1}}$ and $\ang{s_{k+1},s_k}$ 
is crucial; let us assume that this is the first one. Then  
${\rK}={\rQ{s_k,s_{k+1}}}$ is a $1-1$ map of $X_{s_k}$ onto 
$X_{s_{k+1}}$ by \ref{j6}. It follows that the pairs 
$\ang{s_k\we i,s_{k+1}\we i}\;\,(i=0,\,1)$ are crucial and 
the relations ${\rK_i}={\rQ{s_k\we i,\,s_{k+1}\we i}}$ are 
restrictions of $\rK.$ Thus the sets 
$X_{s_{k+1}\we 0}$ and $X_{s_{k+1}\we 1}$ are \dd\rK images 
of resp.\ $X_{s_k\we 0}$ and $X_{s_k\we 1}.$ We now conclude 
that $X_{s_{k+1}\we 0}\cap X_{s_{k+1}\we 1}=\emps$ because 
$X_{s_k\we 0}\cap X_{s_k\we 1}=\emps$ and $\rK$ is $1-1$.


\subsubsection{Reduction to unions of curves}

Let us consider subcase \ref{2c} of Theorem~\ref{m2}. Thus 
$\rGw$ is supposed not to contain pairs $\ang{x,y}$ such that 
$x\in\rL[y]$ or $y\in\rL[x].$ Having this in mind, let us show how 
the construction of Section~\ref{2nd} can be modified to provide 
\ref{x3} in addition to requirements \ref{x1} through \ref{j1}. 

Taking the reasoning in the ``acyclic'' case 
(Subsection~\ref{acy}) as the pattern, we see that the following 
is sufficient. Suppose that $s\in 2^{m-1},\msur$ $u_0=s\we 0,\msur$ 
$u_1=s\we 1.$ Let $u_0=v_0\spp v_1\spp ...\spp v_n=u_1$ be 
the unique sequence of crucial pairs which connects $u_0$  
with $u_1$ in $2^m.$ We have to prove the existence of a system of 
points $x_k\in D_{v_k}$ satisfying $x_0\not=x_n$ and 
$x_k\rQ{k}x_{k+1}$ for all $k<n,$ where 
$\rQ{k}=\rQ{v_kv_{k+1}}$ or resp.\ $\rQ{k}=(\rQ{v_{k+1}v_k})^{-1}$ 
provided $\ang{v_k,v_{k+1}}$ or resp.\ $\ang{v_{k+1},v_k}$ is a 
crucial pair. 

To prove this assertion, it suffices to check the following more 
elementary fact: if $\ang{u,v}$ is a crucial pair and $x\in X_u$ 
(resp.\ $y\in X_v$) then there exist at least two points 
$y\in X_v$ (resp.\ $x\in X_u$) such that 
${x\rQ{uv} y}.$ To prove the fact note that at least one such a 
point $y$ exists by \ref{j4}. Hence if it were unique it would be 
$\od[x],$ therefore $y\in\rL[x],$ contradiction because 
${\rQ{uv}}\sq {\rGw}$.

$\,$\hspep\qed\ (Of Theorem \ref{m2})

\np

\subsection{$\protect\fs12$ graphs in generic model}
\label{s12}

This section presents the proof of Theorem~\ref{tlast}. The 
first subsection shows how to convert \ref2 to an absolute 
statement.  

\subsubsection{Narrow homomorphisms and embeddings}
\label{nar}

Dealing with $\fs12$ graphs, we would like to have the 
absoluteness of the statements: $\Go\hoc\rG$ and $\Go\emc\rG.$ 
This does lot look easily available because the direct 
computation shows that the statements are described by 
$\Sigma^1_4$ formulas, far from the Shonfield level. However, 
two levels can be reduced, at the cost of some specification 
of the nature of homomorphisms and embeddings involved. 

Let us fix a reasonable coding system of countable ordinals 
by reals. Let $\abs z$ denote the ordinal coded by $z,$ and 
$\cO=\ans{z:z\,\hbox{ codes an ordinal}}\,;$ $\cO$ is a 
$\ip11$ set of reals. Let us once and for all associate with 
every $\is12$ graph $\rG$ the representation in the form 
$\rG=\bigcup_{\al<\omi}\rG_\al,$ where the relation 
$x\in\rG_{\abs z}$ is $\ip11$ uniformly on $z\in\cO$ (then the 
components $\rG_\al$ are $\fp11$ of course). In other words 
there exists a $\ip11$ relation $\fpi_{\rG}$ such that 
\dm 
\rG_\al(x,y)\;\llra\;
\fpi_{\rG}(x,y,z)
\hspace{3mm} \hbox{whenever} \hspace{3mm}z\in\cO
\hspace{2mm} \hbox{and} \hspace{2mm}\abs z=\al<\omi\,.
\dm
(Roughly $\fpi_{\rG}(x,y,z)$ says that $z\in\cO$ and we have 
$\psi(x,y,a)$ where $a$ is the \dd{\abs z}th real in the sense 
of the G\"odel wellordering of $\rL[x,y]$ and $\psi(x,y,a)$ is 
a $\Pi^1_1$ formula such that ${\exists\,a\,\psi(x,y,a)}$ 
defines $\rG.$ The Shoenfield absoluteness is used to ensure 
that the whole $\rG$ is captured by the decomposition.) 

Let $\rG$ be a $\is12$ graph. 
A {\it narrow homomorphism\/} of $\Go$ into $\rG$ will be any 
homomorphism (see Introduction) $H$ satisfying the property that, 
for some ordinal $\al_0<\omi,$ we have, for all $a,\,b\in 2^\om,$  
${a\Go b\,\lra\,\ang{H(a),H(b)}\in\bigcup_{\al<\al_0}\rG_\al}.$ 
The notion of a {\it narrow embedding\/} has the same meaning. 

We shall write ${\Go\hocn G}$ (or resp.\ ${\Go\emcn G}$) iff there 
exists a continuous narrow homomorphism (resp.\ embedding) of 
$\Go$ into $\rG$. 

\ble
\label{expl}
Let\/ $\rG$ be a\/ $\is12$ graph on reals. Then\/ 
${\Go\hocn\rG}$ and\/ ${\Go\emcn \rG}$ are\/ $\is12$ relations. 
\ele
\proof Direct compuration. 
\hspep\qed

\noi
The lemma explains the advantage of the ``narrow'' concepts. 
However we have to demonstrate that they are not too narrow 
to violate Theorem~\ref{m2}.

\bpro
\label{enh}
Assume that\/ $\rG$ is a subgraph of a\/ $\is12$ graph\/ $\Gp$ 
in Theorem~\ref{m2}. Then the relations\/ ${\Go\hoc\rG}$ and\/ 
${\Go\emc\rG}$ in\/ {\rmt\ref{mm2}} of Theorem~\ref{m2} can be 
changed to resp.\ ${\Go\hocn\Gp}$ and\/ $\Go\emcn\Gp$.
\epro
\proof Let us come back to Step 1 of the construction described 
in Section~\ref{2nd}. We defined there, for the crucial pair of 
$u_0=\fsg_{m-1}\we 0$ and $u_1=\fsg_{m-1}\we 1$ in $2^m,$ a certain 
nonempty $\od$ set ${\rQ{}}\sq{\rG}.$ Since $\rG\sq\Gp,$ 
there exists an ordinal $\ga_m<\omi$ and a nonempty $\od$ 
set ${\rQ{}'}\sq{\rQ{}}$ such that 
${\rQ{}'}\sq \bigcup_{\al<\ga_m}\Gp_\al.$  
We take now this relation $\rQ{}'$ as the actual $\rQ{}$.  

After the whole construction is finished in this modernized way, 
the ordinal $\al_0=\sup_{m\in\om}\ga_m$ will witness that the 
obtained homomorphism (or embedding, in the modifications used 
in Section~\ref{3rd}) is narrow.\qed

\subsubsection{Extending the universe}
\label{ext}

The proof of Theorem~\ref{tlast} includes the operation of 
embedding of the universe $\rV$ for which the theorem is being 
proved in a larger universe, say $\Vp,$ ``good'' in the sense that 
it satisfies the dichotomy of Theorem~\ref{main}. To facilitate 
this argument, let us pretend that $\rV$ is is fact a countable 
model in a vider $\ZFC$ universe. But this is not a 
serious restrictive assumption --- all that we use it for is the 
existence of certain generic extensions of $\rL,$ the constructible 
subuniverse of $\rV.$ Alternatively, this can be carried out using 
Boolean valued extensions, which does not need any extra hypothesis 
but would reduce the quality of presentation. 

We recall that $\rV,$ the universe for which Theorem~\ref{tlast} 
is being proved, is assumed to be a \dd\zP generic extension of 
$\rL,$ for a forcing notion $\zP\in\rL.$ 
Let $\kpa$ be the cardinal of $\zP$ in $\rL.$ Then 
$\kpa=\aleph_\xi^\rL$ for an ordinal $\xi$.

We let $\Om=\aleph_{\xi+\om}^\rL$ and $\Xi=\aleph_{\xi+\om+1}^\rL.$ 
Obviously $\Om$ and $\Xi$ are cardinals in $\rV$. 

\ble
\label{e}
There is a\/ \dd{\smf\Om}generic extension~\footnote 
{\rmt\ See definitions in Section~\ref{somo}.} 
$\Vp$ of\/ $\rL$ such that\/ $\rV\sq\Vp$.
\ele
\proof Assume, towards contradiction, that such an extension does 
not exist. Then a condition $p\in\zP$ forces over $\rL$ the following: 
$(\ast)\msur$ {\it there does not exist a\/ \dd{\smf\Om}generic 
extension of the universe which is also a\/ \dd{\smf\Om}generic  
extension of\/ $\rL$.}

Let us now consider an arbitrary \dd{\smf\Om}generic extension 
$\Vp$ of $\rL.$ Since $\kpa = (\card\zP)^\rL<\Om,$ there exists a 
\dd\zP generic over $\rL$ set $G\in\Vp$ containing $p.$ We observe 
that $\Vp$ is also an \dd{\smf\Om}generic extension of $\rL[G]$ by 
Proposition~\ref{1sm}, a contradiction with the assumption that 
$p$ forces $(\ast)$.\qed 

\subsubsection{Proof of Theorem \protect\ref{tlast}}
\label{ptlast}
 
Let us fix a $\is12$ graph $\rG,$ defined by a parameter--free 
$\Sigma^1_2$ formula $\vpi(x,y)$ in $\rV.$ 
(As usual, we concentrate on the ``lightface'' case.) 
Avoiding unnecessary complications, 
we let $\Phi(x,y)$ be the $\Sigma^1_2$ formula  
${x\not=y\cj[\hspace{1pt}\vpi(x,y)\orr\vpi(y,x)\hspace{1pt}]\,,}$ 
which defines $\rG$ as well, and explicitly contains the 
requirements of being a graph. $\Gph$ will denote the graph 
${\ans{\ang{x,y}:\Phi(x,y)}}$ in the universe considered  
at the moment; thus for instance $\rG$ is $\Gph$ in $\rV$.

Let $\Vp$ be the extension given by Lemma~\ref e; in particular 
$\Vp$ satisfies $\clps\Om$ and even \ref h of Section~\ref{1st}, 
together with the graph $\Gp=(\Gph)^{\Vp},$ and includes $\rV.$ 
All reals $x\in\rV$ are \dd\Om``weak'' in $\Vp$ by definition.

We now define a subgraph $\rG'\sq\Gp$ to which Theorem~\ref{m2} 
will be applied. 

If $\rG$ {\it is not locally countable\/} in $\rV$ then we 
simply put ${\rG'=\Gp\cap\cWt}$ in $\Vp.$ 

If $\rG$ {\it is locally countable\/} in $\rV$ then we first 
put ${\rG''=\Gp\cap\cWt},$ as above, and then define $\rG'$ 
in $\Vp$ to be the restriction of $\rG''$ on the set $C$ of all 
reals $x\in\cW$ such that the \dd{\Gp}component 
$[x]_{\Gp}=\ans{y:{x \Gp y}}$ is countable in $\Vp$. 

\ble
\label{g'}
$\rG'\sq\cWt\cap\Gp$ is an\/ $\od$ graph in\/ $\Vp.$ 
${\rG={\rG'}\cap\rV.}$ If\/ $\rG$ is acyclic or locally 
countable in\/ $\rV$ then\/ $\rG'$ is acyclic or resp.\ locally 
countable in\/ $\Vp.$
\ele
\proof Let us consider the equality ${\rG={\rG'}\cap\rV.}$ In the 
first version of the definition of $\rG'$ this is a direct 
consequence of the Shoenfield absoluteness theorem. In the locally 
countable version, notice that by the local countability a real 
$x\in\rV$ cannot gain new \dd\Gph neighbours in any extension of 
$\rV$ (here we use the fact that $\Phi$ is $\Sg^1_2$), therefore 
all reals in $\rV$ remain in the domain of $\rG'$. 

The property of being acyclic for a $\is12$ graph can be expressed 
by a $\Pi^1_2$ formula, therefore even $\Gp$ inherits this property 
from $\rG.$ 
(It is not clear whether the local countability can be handled in 
this manner. This is why we changed the definition of $\rG'$ in this 
case so that $\rG'$ acquires this property directly.)\hspep\qed

We have two cases.\vom

{\it Case 2\/}: in $\Vp,$ \underline{not all} reals in $\cW$ 
belong to $\od$ \dd{\rG'}discrete sets. Then $\rG'$ satisfies 
\ref2 of Theorem~\ref{main} in $\Vp$ by Theorem~\ref{m2}. In fact 
we obtain \ref2 in the enhanced form given by 
Proposition~\ref{enh}, that is, ${\Go\hocn\Gp}$ in the main 
part and and ${\Go\emcn\Gp}$ in the ``moreover'' part. Now 
Lemma~\ref{expl} expands this property down to the graph $\rG$ 
in $\rV,$ as required.\vom 

{\it Case 1\/}: in $\Vp,$ \underline{all} reals in $\cW$ 
belong to $\od$ \dd{\rG'}discrete sets. Then $\rG'$ admits an 
$\od$ colouring by countable ordinals in $\Vp$ by Theorem~\ref{m2}. 
Let $\iota(x,\al)$ be an \dd\in formula which defines such a 
colouring in $\Vp.$ It follows from Proposition~\ref{1sm} that, 
since $\Vp$ is a \dd{\smf\Om}generic extension of $\rV,$ the 
restriction $c=c'\res\rV$ is $\od$ in $\rV.$ 
(Indeed each formula with sets in $\rV$ as parameters true in 
$\Vp$ is forced by every forcing condition.) Finally $c$ is an 
$\od$ colouring of $\rG$ by ordinals in $\rV$. 

We have also to prove that $\rG$ admits a $\idh2$ colouring 
in $\rV$ in the case when $\rG$ is locally countable in $\rV.$ This 
is the matter of the next subsection.

\subsubsection{The locally countable case}
\label{lcc}

We keep the notation and setup of Subsection~\ref{ptlast} 
including the assumption of Case 1. We shall also suppose that 
$\rG$ is locally countable in $\rV$.

An elementary reasoning shows that if a locally countable $\od$ 
graph $\rG$ admits an $\od$ colouring by ordinals then it admits 
an $\od$ colouring by {\it countable\/} ordinals. (Indeed, each 
\dd\rG component $[x]_{\rG}$ involves only countably many ordinals 
in the colouring; one simply ``dumps'' all of the ordinals involved 
into $\omi$ in the order preserving way. Notice that the $\is12$ 
graph $x\rG y$ iff $\rL[x]=\rL[y]$ but $x\not=y$ shows that an 
$\od$ colouring by natural numbers may not be available.) 

But this argument does not seem to give a colouring definable 
in $\hc$.

A more elaborate reasoning, based on the next lemma, 
will be used. 

\ble
\label;
The following is true in\/ $\rV.$ 
For any real\/ $x$ there exist$:$ an ordinal\/ $\al'<\omi$ and a 
``term''\/ $\zt'\in\trb{\al'}\cap\rL_{\al'}$ such that the set\/ 
$\zX{\zt'}(\rL_{\al'})$ contains\/ $x$ but does not 
contain any\/ \dd\rG neighbour of\/ $x$.
\ele
Using the least, in the sense of a reasonable wellordering,  
pair $\ang{\al',\zt'}$ satisfying the conclusion of the 
lemma as the colour for $x,$ we get this way a $\idh2$ colouring 
of $\rG$ (in $\rV$) by countable ordinals, as required.\hspep

\proof{}of the lemma. Let us fix a real $x\in\rV.$ In $\Vp,$ there 
exist (see the proof of Theorem~\ref{m1}, item~\ref{mm1}, in 
Section~\ref{1st}): an ordinal $\ga<\Om$ and a ``term'' 
$\zt\in\trb{\ga}\cap\rL$ such that the set $\zX{\zt}=\zX{\zt}(\rL)$ 
(see Subsection~\ref{coding}) contains $x$ and is \dd{\rG'}discrete, 
so does not contain any of \dd{\rG'}neighbours of $x.$ 
Then, since by the local countability all \dd{\Gp}neighbours 
of $x$ in $\Vp$ belong to $\rV$ and are \dd\rG neighbours of 
$x,$ we conclude that $\zX{\zt}$ does not contain any  
\dd\Gp neighbour of $x$ in $\Vp$. 

By the choice of $\Om$ there exists an ordinal $\al<\Om,\msur$ 
$\al>\ga$ such that $\zt\in\rL_\al$ and $\zX{\zt}=\zX{\zt}(\rL_\al).$ 
Since $\Vp$ is a \dd{\smf\Om}generic extension of $\rL[x]$ by 
Proposition~\ref{1sm}, the statement that 
$\zX{\zt}(\rL_\al)$ contains $x$ but does not contain 
any \dd\Gp neighbour of $x$ is \dd{\smf\Om}forced in $\rL[x]$. 

By the local countability again, each \dd\rG 
neighbour of $x$ in $\rV$ belongs to $\rL[x].$  Therefore by the 
L\"ovenheim -- Scolem theorem in $\rL[x]$ there exist ordinals 
$\al'<\Om'<\xi<\omi$ (perhaps $\not<\omi^{\rL[x]}$) and a ``term'' 
$\zt'\in\trb{\al'}\cap\rL_{\al'}$ such that $\rL_\xi[x]$ contains 
all \dd\rG neighbours of $x$ in $\rV,$ models a big enough part 
of $\ZFC,$ and is absolute for $\Phi$ (the $\Sigma^1_2$ formula 
which defines $\rG$ in $\rV$), 
and it is \dd{\smf{\Om'}}forced in $\rL_\xi[x]$ that 
the set $\zX{\zt'}(\rL_{\al'})$ contains $x$ but does not contain 
\dd\Gph neighbours of $x$. 

Let us show that $\al'$ and $\zt'$ satisfy the requirements 
of the lemma. Let $y$ be an arbitrary \dd\rG neighbour of $x$ 
in $\rV;$ we have to prove that $y\not\in\zX{\zt'}(\rL_{\al'})$ 
in $\rV$.

Since $\xi$ is countable in $\rV,$ there exists a 
\dd{\smf{\Om'}}generic extension $\gN$ of $\rL_\xi[x].$ Then it 
is true in $\gN$ that the set 
$\zX{\zt'}(\rL_{\al'})$ contains $x$ but does not contain 
\dd\Gph neighbours of $x$.

Notice that $y\in\gN$ because $\rL_\xi[x]$ 
already contains all \dd\rG neighbours of $x$ in $\rV.$ Furthermore 
it is true in $\gN$ that $y$ is a \dd\Gph neighbour of $x$ by the 
\dd\Phi absoluteness requirement for $\rL_\xi[x].$ Thus 
$y\not\in\zX{\zt'}(\rL_{\al'})$ in $\gN.$ But, since $\al'$ is 
countable in $\gN,$ the statement $y\in\zX{\zt'}(\rL_{\al'})$ is 
actually a $\Sigma^1_1$ formula with some reals in $\gN$ as 
parameters. We conclude that $y\not\in\zX{\zt'}(\rL_{\al'})$ is 
true in $\rV$ as well by the Mostowski absoluteness theorem, as 
required.\hspep\qed

\qed\ (of Theorem~\ref{tlast})

\np

\end{document}